\documentclass[11pt]{article}

\long\def\ig#1{\relax}
\ig{Thanks to Roberto Minio for this def'n.  Compare the def'n of
\comment in AMSTeX.}

\newcount \coefa
\newcount \coefb
\newcount \coefc
\newcount\tempcounta
\newcount\tempcountb
\newcount\tempcountc
\newcount\tempcountd
\newcount\xext
\newcount\yext
\newcount\xoff
\newcount\yoff
\newcount\gap%
\newcount\arrowtypea
\newcount\arrowtypeb
\newcount\arrowtypec
\newcount\arrowtyped
\newcount\arrowtypee
\newcount\height
\newcount\width
\newcount\xpos
\newcount\ypos
\newcount\run
\newcount\rise
\newcount\arrowlength
\newcount\halflength
\newcount\arrowtype
\newdimen\tempdimen
\newdimen\xlen
\newdimen\ylen
\newsavebox{\tempboxa}%
\newsavebox{\tempboxb}%
\newsavebox{\tempboxc}%

\makeatletter
\def\settypes(#1,#2,#3){\arrowtypea#1 \arrowtypeb#2 \arrowtypec#3}
\def\settoheight#1#2{\setbox\@tempboxa\hbox{#2}#1\ht\@tempboxa\relax}%
\def\settodepth#1#2{\setbox\@tempboxa\hbox{#2}#1\dp\@tempboxa\relax}%
\def\settokens[#1`#2`#3`#4]{%
     \def\tokena{#1}\def\tokenb{#2}\def\tokenc{#3}\def\tokend{#4}}
\def\setsqparms[#1`#2`#3`#4;#5`#6]{%
\arrowtypea #1
\arrowtypeb #2
\arrowtypec #3
\arrowtyped #4
\width #5
\height #6
}
\def\setpos(#1,#2){\xpos=#1 \ypos#2}

\def\bfig{\begin{picture}(\xext,\yext)(\xoff,\yoff)}
\def\efig{\end{picture}}

\def\putbox(#1,#2)#3{\put(#1,#2){\makebox(0,0){$#3$}}}

\def\settriparms[#1`#2`#3;#4]{\settripairparms[#1`#2`#3`1`1;#4]}%

\def\settripairparms[#1`#2`#3`#4`#5;#6]{%
\arrowtypea #1
\arrowtypeb #2
\arrowtypec #3
\arrowtyped #4
\arrowtypee #5
\width #6
\height #6
}

\def\resetparms{\settripairparms[1`1`1`1`1;500]\width 500}

\resetparms

\def\mvector(#1,#2)#3{
\put(0,0){\vector(#1,#2){#3}}%
\put(0,0){\vector(#1,#2){30}}%
}
\def\evector(#1,#2)#3{{
\arrowlength #3
\put(0,0){\vector(#1,#2){\arrowlength}}%
\advance \arrowlength by-30
\put(0,0){\vector(#1,#2){\arrowlength}}%
}}

\def\horsize#1#2{%
\settowidth{\tempdimen}{$#2$}%
#1=\tempdimen
\divide #1 by\unitlength
}

\def\vertsize#1#2{%
\settoheight{\tempdimen}{$#2$}%
#1=\tempdimen
\settodepth{\tempdimen}{$#2$}%
\advance #1 by\tempdimen
\divide #1 by\unitlength
}

\def\vertadjust[#1`#2`#3]{%
\vertsize{\tempcounta}{#1}%
\vertsize{\tempcountb}{#2}%
\ifnum \tempcounta<\tempcountb \tempcounta=\tempcountb \fi
\divide\tempcounta by2
\vertsize{\tempcountb}{#3}%
\ifnum \tempcountb>0 \advance \tempcountb by20 \fi
\ifnum \tempcounta<\tempcountb \tempcounta=\tempcountb \fi
}

\def\horadjust[#1`#2`#3]{%
\horsize{\tempcounta}{#1}%
\horsize{\tempcountb}{#2}%
\ifnum \tempcounta<\tempcountb \tempcounta=\tempcountb \fi
\divide\tempcounta by20
\horsize{\tempcountb}{#3}%
\ifnum \tempcountb>0 \advance \tempcountb by60 \fi
\ifnum \tempcounta<\tempcountb \tempcounta=\tempcountb \fi
}

\ig{ In this procedure, #1 is the paramater that sticks out all the way,
#2 sticks out the least and #3 is a label sticking out half way.  #4 is
the amount of the offset.}

\def\sladjust[#1`#2`#3]#4{%
\tempcountc=#4
\horsize{\tempcounta}{#1}%
\divide \tempcounta by2
\horsize{\tempcountb}{#2}%
\divide \tempcountb by2
\advance \tempcountb by-\tempcountc
\ifnum \tempcounta<\tempcountb \tempcounta=\tempcountb\fi
\divide \tempcountc by2
\horsize{\tempcountb}{#3}%
\advance \tempcountb by-\tempcountc
\ifnum \tempcountb>0 \advance \tempcountb by80\fi
\ifnum \tempcounta<\tempcountb \tempcounta=\tempcountb\fi
\advance\tempcounta by20
}

\def\putvector(#1,#2)(#3,#4)#5#6{{%
\xpos=#1
\ypos=#2
\run=#3
\rise=#4
\arrowlength=#5
\arrowtype=#6
\ifnum \arrowtype<0
    \ifnum \run=0
        \advance \ypos by-\arrowlength
    \else
        \tempcounta \arrowlength
        \multiply \tempcounta by\rise
        \divide \tempcounta by\run
        \ifnum\run>0
            \advance \xpos by\arrowlength
            \advance \ypos by\tempcounta
        \else
            \advance \xpos by-\arrowlength
            \advance \ypos by-\tempcounta
        \fi
    \fi
    \multiply \arrowtype by-1
    \multiply \rise by-1
    \multiply \run by-1
\fi
\ifnum \arrowtype=1
    \put(\xpos,\ypos){\vector(\run,\rise){\arrowlength}}%
\else\ifnum \arrowtype=2
    \put(\xpos,\ypos){\mvector(\run,\rise)\arrowlength}%
\else\ifnum\arrowtype=3
    \put(\xpos,\ypos){\evector(\run,\rise){\arrowlength}}%
\fi\fi\fi
}}

\def\putsplitvector(#1,#2)#3#4{
\xpos #1
\ypos #2
\arrowtype #4
\halflength #3
\arrowlength #3
\gap 140
\advance \halflength by-\gap
\divide \halflength by2
\ifnum \arrowtype=1
    \put(\xpos,\ypos){\line(0,-1){\halflength}}%
    \advance\ypos by-\halflength
    \advance\ypos by-\gap
    \put(\xpos,\ypos){\vector(0,-1){\halflength}}%
\else\ifnum \arrowtype=2
    \put(\xpos,\ypos){\line(0,-1)\halflength}%
    \put(\xpos,\ypos){\vector(0,-1)3}%
    \advance\ypos by-\halflength
    \advance\ypos by-\gap
    \put(\xpos,\ypos){\vector(0,-1){\halflength}}%
\else\ifnum\arrowtype=3
    \put(\xpos,\ypos){\line(0,-1)\halflength}%
    \advance\ypos by-\halflength
    \advance\ypos by-\gap
    \put(\xpos,\ypos){\evector(0,-1){\halflength}}%
\else\ifnum \arrowtype=-1
    \advance \ypos by-\arrowlength
    \put(\xpos,\ypos){\line(0,1){\halflength}}%
    \advance\ypos by\halflength
    \advance\ypos by\gap
    \put(\xpos,\ypos){\vector(0,1){\halflength}}%
\else\ifnum \arrowtype=-2
    \advance \ypos by-\arrowlength
    \put(\xpos,\ypos){\line(0,1)\halflength}%
    \put(\xpos,\ypos){\vector(0,1)3}%
    \advance\ypos by\halflength
    \advance\ypos by\gap
    \put(\xpos,\ypos){\vector(0,1){\halflength}}%
\else\ifnum\arrowtype=-3
    \advance \ypos by-\arrowlength
    \put(\xpos,\ypos){\line(0,1)\halflength}%
    \advance\ypos by\halflength
    \advance\ypos by\gap
    \put(\xpos,\ypos){\evector(0,1){\halflength}}%
\fi\fi\fi\fi\fi\fi
}

\def\putmorphism(#1)(#2,#3)[#4`#5`#6]#7#8#9{{%
\run #2
\rise #3
\ifnum\rise=0
  \puthmorphism(#1)[#4`#5`#6]{#7}{#8}{#9}%
\else\ifnum\run=0
  \putvmorphism(#1)[#4`#5`#6]{#7}{#8}{#9}%
\else
\setpos(#1)%
\arrowlength #7
\arrowtype #8
\ifnum\run=0
\else\ifnum\rise=0
\else
\ifnum\run>0
    \coefa=1
\else
   \coefa=-1
\fi
\ifnum\arrowtype>0
   \coefb=0
   \coefc=-1
\else
   \coefb=\coefa
   \coefc=1
   \arrowtype=-\arrowtype
\fi
\width=2
\multiply \width by\run
\divide \width by\rise
\ifnum \width<0  \width=-\width\fi
\advance\width by60
\if l#9 \width=-\width\fi
\putbox(\xpos,\ypos){#4}
{\multiply \coefa by\arrowlength
\advance\xpos by\coefa
\multiply \coefa by\rise
\divide \coefa by\run
\advance \ypos by\coefa
\putbox(\xpos,\ypos){#5} }%
{\multiply \coefa by\arrowlength
\divide \coefa by2
\advance \xpos by\coefa
\advance \xpos by\width
\multiply \coefa by\rise
\divide \coefa by\run
\advance \ypos by\coefa
\if l#9%
   \put(\xpos,\ypos){\makebox(0,0)[r]{$#6$}}%
\else\if r#9%
   \put(\xpos,\ypos){\makebox(0,0)[l]{$#6$}}%
\fi\fi }%
{\multiply \rise by-\coefc
\multiply \run by-\coefc
\multiply \coefb by\arrowlength
\advance \xpos by\coefb
\multiply \coefb by\rise
\divide \coefb by\run
\advance \ypos by\coefb
\multiply \coefc by70
\advance \ypos by\coefc
\multiply \coefc by\run
\divide \coefc by\rise
\advance \xpos by\coefc
\multiply \coefa by140
\multiply \coefa by\run
\divide \coefa by\rise
\advance \arrowlength by\coefa
\ifnum \arrowtype=1
   \put(\xpos,\ypos){\vector(\run,\rise){\arrowlength}}%
\else\ifnum\arrowtype=2
   \put(\xpos,\ypos){\mvector(\run,\rise){\arrowlength}}%
\else\ifnum\arrowtype=3
   \put(\xpos,\ypos){\evector(\run,\rise){\arrowlength}}%
\fi\fi\fi}\fi\fi\fi\fi}}

\def\puthmorphism(#1,#2)[#3`#4`#5]#6#7#8{{%
\xpos #1
\ypos #2
\width #6
\arrowlength #6
\putbox(\xpos,\ypos){#3\vphantom{#4}}%
{\advance \xpos by\arrowlength
\putbox(\xpos,\ypos){\vphantom{#3}#4}}%
\horsize{\tempcounta}{#3}%
\horsize{\tempcountb}{#4}%
\divide \tempcounta by2
\divide \tempcountb by2
\advance \tempcounta by30
\advance \tempcountb by30
\advance \xpos by\tempcounta
\advance \arrowlength by-\tempcounta
\advance \arrowlength by-\tempcountb
\putvector(\xpos,\ypos)(1,0){\arrowlength}{#7}%
\divide \arrowlength by2
\advance \xpos by\arrowlength
\vertsize{\tempcounta}{#5}%
\divide\tempcounta by2
\advance \tempcounta by20
\if a#8 %
   \advance \ypos by\tempcounta
   \putbox(\xpos,\ypos){#5}%
\else
   \advance \ypos by-\tempcounta
   \putbox(\xpos,\ypos){#5}%
\fi}}

\def\putvmorphism(#1,#2)[#3`#4`#5]#6#7#8{{%
\xpos #1
\ypos #2
\arrowlength #6
\arrowtype #7
\settowidth{\xlen}{$#5$}%
\putbox(\xpos,\ypos){#3}%
{\advance \ypos by-\arrowlength
\putbox(\xpos,\ypos){#4}}%
{\advance\arrowlength by-140
\advance \ypos by-70
\ifdim\xlen>0pt
   \if m#8%
      \putsplitvector(\xpos,\ypos){\arrowlength}{\arrowtype}%
   \else
      \putvector(\xpos,\ypos)(0,-1){\arrowlength}{\arrowtype}%
   \fi
\else
   \putvector(\xpos,\ypos)(0,-1){\arrowlength}{\arrowtype}%
\fi}%
\ifdim\xlen>0pt
   \divide \arrowlength by2
   \advance\ypos by-\arrowlength
   \if l#8%
      \advance \xpos by-40
      \put(\xpos,\ypos){\makebox(0,0)[r]{$#5$}}%
   \else\if r#8%
      \advance \xpos by40
      \put(\xpos,\ypos){\makebox(0,0)[l]{$#5$}}%
   \else
      \putbox(\xpos,\ypos){#5}%
   \fi\fi
\fi
}}

\def\topadjust[#1`#2`#3]{%
\yoff=10
\vertadjust[#1`#2`{#3}]%
\advance \yext by\tempcounta
\advance \yext by 10
}
\def\botadjust[#1`#2`#3]{%
\vertadjust[#1`#2`{#3}]%
\advance \yext by\tempcounta
\advance \yoff by-\tempcounta
}
\def\leftadjust[#1`#2`#3]{%
\xoff=0
\horadjust[#1`#2`{#3}]%
\advance \xext by\tempcounta
\advance \xoff by-\tempcounta
}
\def\rightadjust[#1`#2`#3]{%
\horadjust[#1`#2`{#3}]%
\advance \xext by\tempcounta
}
\def\rightsladjust[#1`#2`#3]{%
\sladjust[#1`#2`{#3}]{\width}%
\advance \xext by\tempcounta
}
\def\leftsladjust[#1`#2`#3]{%
\xoff=0
\sladjust[#1`#2`{#3}]{\width}%
\advance \xext by\tempcounta
\advance \xoff by-\tempcounta
}
\def\adjust[#1`#2;#3`#4;#5`#6;#7`#8]{%
\topadjust[#1``{#2}]
\leftadjust[#3``{#4}]
\rightadjust[#5``{#6}]
\botadjust[#7``{#8}]}

\def\putsquarep<#1>(#2)[#3;#4`#5`#6`#7]{{%
\setsqparms[#1]%
\setpos(#2)%
\settokens[#3]%
\puthmorphism(\xpos,\ypos)[\tokenc`\tokend`{#7}]{\width}{\arrowtyped}b%
\advance\ypos by \height
\puthmorphism(\xpos,\ypos)[\tokena`\tokenb`{#4}]{\width}{\arrowtypea}a%
\putvmorphism(\xpos,\ypos)[``{#5}]{\height}{\arrowtypeb}l%
\advance\xpos by \width
\putvmorphism(\xpos,\ypos)[``{#6}]{\height}{\arrowtypec}r%
}}

\def\putsquare{\@ifnextchar <{\putsquarep}{\putsquarep%
   <\arrowtypea`\arrowtypeb`\arrowtypec`\arrowtyped;\width`\height>}}
\def\square{\@ifnextchar< {\squarep}{\squarep
   <\arrowtypea`\arrowtypeb`\arrowtypec`\arrowtyped;\width`\height>}}
\def\squarep<#1>[#2`#3`#4`#5;#6`#7`#8`#9]{{
\setsqparms[#1]
\xext=\width                                          
\yext=\height                                         
\topadjust[#2`#3`{#6}]
\botadjust[#4`#5`{#9}]
\leftadjust[#2`#4`{#7}]
\rightadjust[#3`#5`{#8}]
\begin{picture}(\xext,\yext)(\xoff,\yoff)
\putsquarep<\arrowtypea`\arrowtypeb`\arrowtypec`\arrowtyped;\width`\height>%
(0,0)[#2`#3`#4`#5;#6`#7`#8`{#9}]%
\end{picture}%
}}

\def\putptrianglep<#1>(#2,#3)[#4`#5`#6;#7`#8`#9]{{%
\settriparms[#1]%
\xpos=#2 \ypos=#3
\advance\ypos by \height
\puthmorphism(\xpos,\ypos)[#4`#5`{#7}]{\height}{\arrowtypea}a%
\putvmorphism(\xpos,\ypos)[`#6`{#8}]{\height}{\arrowtypeb}l%
\advance\xpos by\height
\putmorphism(\xpos,\ypos)(-1,-1)[``{#9}]{\height}{\arrowtypec}r%
}}

\def\putptriangle{\@ifnextchar <{\putptrianglep}{\putptrianglep
   <\arrowtypea`\arrowtypeb`\arrowtypec;\height>}}
\def\ptriangle{\@ifnextchar <{\ptrianglep}{\ptrianglep
   <\arrowtypea`\arrowtypeb`\arrowtypec;\height>}}

\def\ptrianglep<#1>[#2`#3`#4;#5`#6`#7]{{
\settriparms[#1]%
\width=\height                         
\xext=\width                           
\yext=\width                           
\topadjust[#2`#3`{#5}]
\botadjust[#3``]
\leftadjust[#2`#4`{#6}]
\rightsladjust[#3`#4`{#7}]
\begin{picture}(\xext,\yext)(\xoff,\yoff)
\putptrianglep<\arrowtypea`\arrowtypeb`\arrowtypec;\height>%
(0,0)[#2`#3`#4;#5`#6`{#7}]%
\end{picture}%
}}

\def\putqtrianglep<#1>(#2,#3)[#4`#5`#6;#7`#8`#9]{{%
\settriparms[#1]%
\xpos=#2 \ypos=#3
\advance\ypos by\height
\puthmorphism(\xpos,\ypos)[#4`#5`{#7}]{\height}{\arrowtypea}a%
\putmorphism(\xpos,\ypos)(1,-1)[``{#8}]{\height}{\arrowtypeb}l%
\advance\xpos by\height
\putvmorphism(\xpos,\ypos)[`#6`{#9}]{\height}{\arrowtypec}r%
}}

\def\putqtriangle{\@ifnextchar <{\putqtrianglep}{\putqtrianglep
   <\arrowtypea`\arrowtypeb`\arrowtypec;\height>}}
\def\qtriangle{\@ifnextchar <{\qtrianglep}{\qtrianglep
   <\arrowtypea`\arrowtypeb`\arrowtypec;\height>}}

\def\qtrianglep<#1>[#2`#3`#4;#5`#6`#7]{{
\settriparms[#1]
\width=\height                         
\xext=\width                           
\yext=\height                          
\topadjust[#2`#3`{#5}]
\botadjust[#4``]
\leftsladjust[#2`#4`{#6}]
\rightadjust[#3`#4`{#7}]
\begin{picture}(\xext,\yext)(\xoff,\yoff)
\putqtrianglep<\arrowtypea`\arrowtypeb`\arrowtypec;\height>%
(0,0)[#2`#3`#4;#5`#6`{#7}]%
\end{picture}%
}}

\def\putdtrianglep<#1>(#2,#3)[#4`#5`#6;#7`#8`#9]{{%
\settriparms[#1]%
\xpos=#2 \ypos=#3
\puthmorphism(\xpos,\ypos)[#5`#6`{#9}]{\height}{\arrowtypec}b%
\advance\xpos by \height \advance\ypos by\height
\putmorphism(\xpos,\ypos)(-1,-1)[``{#7}]{\height}{\arrowtypea}l%
\putvmorphism(\xpos,\ypos)[#4``{#8}]{\height}{\arrowtypeb}r%
}}

\def\putdtriangle{\@ifnextchar <{\putdtrianglep}{\putdtrianglep
   <\arrowtypea`\arrowtypeb`\arrowtypec;\height>}}
\def\dtriangle{\@ifnextchar <{\dtrianglep}{\dtrianglep
   <\arrowtypea`\arrowtypeb`\arrowtypec;\height>}}

\def\dtrianglep<#1>[#2`#3`#4;#5`#6`#7]{{
\settriparms[#1]
\width=\height                         
\xext=\width                           
\yext=\height                          
\topadjust[#2``]
\botadjust[#3`#4`{#7}]
\leftsladjust[#3`#2`{#5}]
\rightadjust[#2`#4`{#6}]
\begin{picture}(\xext,\yext)(\xoff,\yoff)
\putdtrianglep<\arrowtypea`\arrowtypeb`\arrowtypec;\height>%
(0,0)[#2`#3`#4;#5`#6`{#7}]%
\end{picture}%
}}

\def\putbtrianglep<#1>(#2,#3)[#4`#5`#6;#7`#8`#9]{{%
\settriparms[#1]%
\xpos=#2 \ypos=#3
\puthmorphism(\xpos,\ypos)[#5`#6`{#9}]{\height}{\arrowtypec}b%
\advance\ypos by\height
\putmorphism(\xpos,\ypos)(1,-1)[``{#8}]{\height}{\arrowtypeb}r%
\putvmorphism(\xpos,\ypos)[#4``{#7}]{\height}{\arrowtypea}l%
}}

\def\putbtriangle{\@ifnextchar <{\putbtrianglep}{\putbtrianglep
   <\arrowtypea`\arrowtypeb`\arrowtypec;\height>}}
\def\btriangle{\@ifnextchar <{\btrianglep}{\btrianglep
   <\arrowtypea`\arrowtypeb`\arrowtypec;\height>}}

\def\btrianglep<#1>[#2`#3`#4;#5`#6`#7]{{
\settriparms[#1]
\width=\height                         
\xext=\width                           
\yext=\height                          
\topadjust[#2``]
\botadjust[#3`#4`{#7}]
\leftadjust[#2`#3`{#5}]
\rightsladjust[#4`#2`{#6}]
\begin{picture}(\xext,\yext)(\xoff,\yoff)
\putbtrianglep<\arrowtypea`\arrowtypeb`\arrowtypec;\height>%
(0,0)[#2`#3`#4;#5`#6`{#7}]%
\end{picture}%
}}

\def\putAtrianglep<#1>(#2,#3)[#4`#5`#6;#7`#8`#9]{{%
\settriparms[#1]%
\xpos=#2 \ypos=#3
{\multiply \height by2
\puthmorphism(\xpos,\ypos)[#5`#6`{#9}]{\height}{\arrowtypec}b}%
\advance\xpos by\height \advance\ypos by\height
\putmorphism(\xpos,\ypos)(-1,-1)[#4``{#7}]{\height}{\arrowtypea}l%
\putmorphism(\xpos,\ypos)(1,-1)[``{#8}]{\height}{\arrowtypeb}r%
}}

\def\putAtriangle{\@ifnextchar <{\putAtrianglep}{\putAtrianglep
   <\arrowtypea`\arrowtypeb`\arrowtypec;\height>}}
\def\Atriangle{\@ifnextchar <{\Atrianglep}{\Atrianglep
   <\arrowtypea`\arrowtypeb`\arrowtypec;\height>}}

\def\Atrianglep<#1>[#2`#3`#4;#5`#6`#7]{{
\settriparms[#1]
\width=\height                         
\xext=\width                           
\yext=\height                          
\topadjust[#2``]
\botadjust[#3`#4`{#7}]
\multiply \xext by2 
\leftsladjust[#3`#2`{#5}]
\rightsladjust[#4`#2`{#6}]
\begin{picture}(\xext,\yext)(\xoff,\yoff)%
\putAtrianglep<\arrowtypea`\arrowtypeb`\arrowtypec;\height>%
(0,0)[#2`#3`#4;#5`#6`{#7}]%
\end{picture}%
}}

\def\putAtrianglepairp<#1>(#2)[#3;#4`#5`#6`#7`#8]{{
\settripairparms[#1]%
\setpos(#2)%
\settokens[#3]%
\puthmorphism(\xpos,\ypos)[\tokenb`\tokenc`{#7}]{\height}{\arrowtyped}b%
\advance\xpos by\height
\advance\ypos by\height
\putmorphism(\xpos,\ypos)(-1,-1)[\tokena``{#4}]{\height}{\arrowtypea}l%
\putvmorphism(\xpos,\ypos)[``{#5}]{\height}{\arrowtypeb}m%
\putmorphism(\xpos,\ypos)(1,-1)[``{#6}]{\height}{\arrowtypec}r%
}}

\def\putAtrianglepair{\@ifnextchar <{\putAtrianglepairp}{\putAtrianglepairp%
   <\arrowtypea`\arrowtypeb`\arrowtypec`\arrowtyped`\arrowtypee;\height>}}
\def\Atrianglepair{\@ifnextchar <{\Atrianglepairp}{\Atrianglepairp%
   <\arrowtypea`\arrowtypeb`\arrowtypec`\arrowtyped`\arrowtypee;\height>}}

\def\Atrianglepairp<#1>[#2;#3`#4`#5`#6`#7]{{%
\settripairparms[#1]%
\settokens[#2]%
\width=\height
\xext=\width
\yext=\height
\topadjust[\tokena``]%
\vertadjust[\tokenb`\tokenc`{#6}]
\tempcountd=\tempcounta                       
\vertadjust[\tokenc`\tokend`{#7}]
\ifnum\tempcounta<\tempcountd                 
\tempcounta=\tempcountd\fi                    
\advance \yext by\tempcounta                  
\advance \yoff by-\tempcounta                 %
\multiply \xext by2 
\leftsladjust[\tokenb`\tokena`{#3}]
\rightsladjust[\tokend`\tokena`{#5}]%
\begin{picture}(\xext,\yext)(\xoff,\yoff)%
\putAtrianglepairp
<\arrowtypea`\arrowtypeb`\arrowtypec`\arrowtyped`\arrowtypee;\height>%
(0,0)[#2;#3`#4`#5`#6`{#7}]%
\end{picture}%
}}

\def\putVtrianglep<#1>(#2,#3)[#4`#5`#6;#7`#8`#9]{{%
\settriparms[#1]%
\xpos=#2 \ypos=#3
\advance\ypos by\height
{\multiply\height by2
\puthmorphism(\xpos,\ypos)[#4`#5`{#7}]{\height}{\arrowtypea}a}%
\putmorphism(\xpos,\ypos)(1,-1)[`#6`{#8}]{\height}{\arrowtypeb}l%
\advance\xpos by\height
\advance\xpos by\height
\putmorphism(\xpos,\ypos)(-1,-1)[``{#9}]{\height}{\arrowtypec}r%
}}

\def\putVtriangle{\@ifnextchar <{\putVtrianglep}{\putVtrianglep
   <\arrowtypea`\arrowtypeb`\arrowtypec;\height>}}
\def\Vtriangle{\@ifnextchar <{\Vtrianglep}{\Vtrianglep
   <\arrowtypea`\arrowtypeb`\arrowtypec;\height>}}

\def\Vtrianglep<#1>[#2`#3`#4;#5`#6`#7]{{
\settriparms[#1]
\width=\height                         
\xext=\width                           
\yext=\height                          
\topadjust[#2`#3`{#5}]
\botadjust[#4``]
\multiply \xext by2 
\leftsladjust[#2`#3`{#6}]
\rightsladjust[#3`#4`{#7}]
\begin{picture}(\xext,\yext)(\xoff,\yoff)%
\putVtrianglep<\arrowtypea`\arrowtypeb`\arrowtypec;\height>%
(0,0)[#2`#3`#4;#5`#6`{#7}]%
\end{picture}%
}}

\def\putVtrianglepairp<#1>(#2)[#3;#4`#5`#6`#7`#8]{{
\settripairparms[#1]%
\setpos(#2)%
\settokens[#3]%
\advance\ypos by\height
\putmorphism(\xpos,\ypos)(1,-1)[`\tokend`{#6}]{\height}{\arrowtypec}l%
\puthmorphism(\xpos,\ypos)[\tokena`\tokenb`{#4}]{\height}{\arrowtypea}a%
\advance\xpos by\height
\putvmorphism(\xpos,\ypos)[``{#7}]{\height}{\arrowtyped}m%
\advance\xpos by\height
\putmorphism(\xpos,\ypos)(-1,-1)[``{#8}]{\height}{\arrowtypee}r%
}}

\def\putVtrianglepair{\@ifnextchar <{\putVtrianglepairp}{\putVtrianglepairp%
    <\arrowtypea`\arrowtypeb`\arrowtypec`\arrowtyped`\arrowtypee;\height>}}
\def\Vtrianglepair{\@ifnextchar <{\Vtrianglepairp}{\Vtrianglepairp%
    <\arrowtypea`\arrowtypeb`\arrowtypec`\arrowtyped`\arrowtypee;\height>}}

\def\Vtrianglepairp<#1>[#2;#3`#4`#5`#6`#7]{{%
\settripairparms[#1]%
\settokens[#2]
\xext=\height                  
\width=\height                 
\yext=\height                  
\vertadjust[\tokena`\tokenb`{#4}]
\tempcountd=\tempcounta        
\vertadjust[\tokenb`\tokenc`{#5}]
\ifnum\tempcounta<\tempcountd%
\tempcounta=\tempcountd\fi
\advance \yext by\tempcounta
\botadjust[\tokend``]%
\multiply \xext by2
\leftsladjust[\tokena`\tokend`{#6}]%
\rightsladjust[\tokenc`\tokend`{#7}]%
\begin{picture}(\xext,\yext)(\xoff,\yoff)%
\putVtrianglepairp
<\arrowtypea`\arrowtypeb`\arrowtypec`\arrowtyped`\arrowtypee;\height>%
(0,0)[#2;#3`#4`#5`#6`{#7}]%
\end{picture}%
}}

\def\putCtrianglep<#1>(#2,#3)[#4`#5`#6;#7`#8`#9]{{%
\settriparms[#1]%
\xpos=#2 \ypos=#3
\advance\ypos by\height
\putmorphism(\xpos,\ypos)(1,-1)[``{#9}]{\height}{\arrowtypec}l%
\advance\xpos by\height
\advance\ypos by\height
\putmorphism(\xpos,\ypos)(-1,-1)[#4`#5`{#7}]{\height}{\arrowtypea}l%
{\multiply\height by 2
\putvmorphism(\xpos,\ypos)[`#6`{#8}]{\height}{\arrowtypeb}r}%
}}

\def\putCtriangle{\@ifnextchar <{\putCtrianglep}{\putCtrianglep
    <\arrowtypea`\arrowtypeb`\arrowtypec;\height>}}
\def\Ctriangle{\@ifnextchar <{\Ctrianglep}{\Ctrianglep
    <\arrowtypea`\arrowtypeb`\arrowtypec;\height>}}

\def\Ctrianglep<#1>[#2`#3`#4;#5`#6`#7]{{
\settriparms[#1]
\width=\height                          
\xext=\width                            
\yext=\height                           
\multiply \yext by2 
\topadjust[#2``]
\botadjust[#4``]
\sladjust[#3`#2`{#5}]{\width}
\tempcountd=\tempcounta                 
\sladjust[#3`#4`{#7}]{\width}
\ifnum \tempcounta<\tempcountd          
\tempcounta=\tempcountd\fi              
\advance \xext by\tempcounta            
\advance \xoff by-\tempcounta           %
\rightadjust[#2`#4`{#6}]
\begin{picture}(\xext,\yext)(\xoff,\yoff)%
\putCtrianglep<\arrowtypea`\arrowtypeb`\arrowtypec;\height>%
(0,0)[#2`#3`#4;#5`#6`{#7}]%
\end{picture}%
}}

\def\putDtrianglep<#1>(#2,#3)[#4`#5`#6;#7`#8`#9]{{%
\settriparms[#1]%
\xpos=#2 \ypos=#3
\advance\xpos by\height \advance\ypos by\height
\putmorphism(\xpos,\ypos)(-1,-1)[``{#9}]{\height}{\arrowtypec}r%
\advance\xpos by-\height \advance\ypos by\height
\putmorphism(\xpos,\ypos)(1,-1)[`#5`{#8}]{\height}{\arrowtypeb}r%
{\multiply\height by 2
\putvmorphism(\xpos,\ypos)[#4`#6`{#7}]{\height}{\arrowtypea}l}%
}}

\def\putDtriangle{\@ifnextchar <{\putDtrianglep}{\putDtrianglep
    <\arrowtypea`\arrowtypeb`\arrowtypec;\height>}}
\def\Dtriangle{\@ifnextchar <{\Dtrianglep}{\Dtrianglep
   <\arrowtypea`\arrowtypeb`\arrowtypec;\height>}}

\def\Dtrianglep<#1>[#2`#3`#4;#5`#6`#7]{{
\settriparms[#1]
\width=\height                         
\xext=\height                          
\yext=\height                          
\multiply \yext by2 
\topadjust[#2``]
\botadjust[#4``]
\leftadjust[#2`#4`{#5}]
\sladjust[#3`#2`{#5}]{\height}
\tempcountd=\tempcountd                
\sladjust[#3`#4`{#7}]{\height}
\ifnum \tempcounta<\tempcountd         
\tempcounta=\tempcountd\fi             
\advance \xext by\tempcounta           %
\begin{picture}(\xext,\yext)(\xoff,\yoff)
\putDtrianglep<\arrowtypea`\arrowtypeb`\arrowtypec;\height>%
(0,0)[#2`#3`#4;#5`#6`{#7}]%
\end{picture}%
}}

\def\setrecparms[#1`#2]{\width=#1 \height=#2}%
\def\recursep<#1`#2>[#3;#4`#5`#6`#7`#8]{{%
\width=#1 \height=#2
\settokens[#3]
\settowidth{\tempdimen}{$\tokena$}
\ifdim\tempdimen=0pt
  \savebox{\tempboxa}{\hbox{$\tokenb$}}%
  \savebox{\tempboxb}{\hbox{$\tokend$}}%
  \savebox{\tempboxc}{\hbox{$#6$}}%
\else
  \savebox{\tempboxa}{\hbox{$\hbox{$\tokena$}\times\hbox{$\tokenb$}$}}%
  \savebox{\tempboxb}{\hbox{$\hbox{$\tokena$}\times\hbox{$\tokend$}$}}%
  \savebox{\tempboxc}{\hbox{$\hbox{$\tokena$}\times\hbox{$#6$}$}}%
\fi
\ypos=\height
\divide\ypos by 2
\xpos=\ypos
\advance\xpos by \width
\xext=\xpos \yext=\height
\topadjust[#3`\usebox{\tempboxa}`{#4}]%
\botadjust[#5`\usebox{\tempboxb}`{#8}]%
\sladjust[\tokenc`\tokenb`{#5}]{\ypos}%
\tempcountd=\tempcounta
\sladjust[\tokenc`\tokend`{#5}]{\ypos}%
\ifnum \tempcounta<\tempcountd
\tempcounta=\tempcountd\fi
\advance \xext by\tempcounta
\advance \xoff by-\tempcounta
\rightadjust[\usebox{\tempboxa}`\usebox{\tempboxb}`\usebox{\tempboxc}]%
\bfig
\putCtrianglep<-1`1`1;\ypos>(0,0)[`\tokenc`;#5`#6`{#7}]%
\puthmorphism(\ypos,0)[\tokend`\usebox{\tempboxb}`{#8}]{\width}{-1}b%
\puthmorphism(\ypos,\height)[\tokenb`\usebox{\tempboxa}`{#4}]{\width}{-1}a%
\advance\ypos by \width
\putvmorphism(\ypos,\height)[``\usebox{\tempboxc}]{\height}1r%
\efig
}}

\def\recurse{\@ifnextchar <{\recursep}{\recursep<\width`\height>}}

\def\puttwohmorphisms(#1,#2)[#3`#4;#5`#6]#7#8#9{{%
\puthmorphism(#1,#2)[#3`#4`]{#7}0a
\ypos=#2
\advance\ypos by 20
\puthmorphism(#1,\ypos)[\phantom{#3}`\phantom{#4}`#5]{#7}{#8}a
\advance\ypos by -40
\puthmorphism(#1,\ypos)[\phantom{#3}`\phantom{#4}`#6]{#7}{#9}b
}}

\def\puttwovmorphisms(#1,#2)[#3`#4;#5`#6]#7#8#9{{%
\putvmorphism(#1,#2)[#3`#4`]{#7}0a
\xpos=#1
\advance\xpos by -20
\putvmorphism(\xpos,#2)[\phantom{#3}`\phantom{#4}`#5]{#7}{#8}l
\advance\xpos by 40
\putvmorphism(\xpos,#2)[\phantom{#3}`\phantom{#4}`#6]{#7}{#9}r
}}

\def\puthcoequalizer(#1)[#2`#3`#4;#5`#6`#7]#8#9{{%
\setpos(#1)%
\puttwohmorphisms(\xpos,\ypos)[#2`#3;#5`#6]{#8}11%
\advance\xpos by #8
\puthmorphism(\xpos,\ypos)[\phantom{#3}`#4`#7]{#8}1{#9}
}}

\def\putvcoequalizer(#1)[#2`#3`#4;#5`#6`#7]#8#9{{%
\setpos(#1)%
\puttwovmorphisms(\xpos,\ypos)[#2`#3;#5`#6]{#8}11%
\advance\ypos by -#8
\putvmorphism(\xpos,\ypos)[\phantom{#3}`#4`#7]{#8}1{#9}
}}

\def\putthreehmorphisms(#1)[#2`#3;#4`#5`#6]#7(#8)#9{{%
\setpos(#1) \settypes(#8)
\if a#9 %
     \vertsize{\tempcounta}{#5}%
     \vertsize{\tempcountb}{#6}%
     \ifnum \tempcounta<\tempcountb \tempcounta=\tempcountb \fi
\else
     \vertsize{\tempcounta}{#4}%
     \vertsize{\tempcountb}{#5}%
     \ifnum \tempcounta<\tempcountb \tempcounta=\tempcountb \fi
\fi
\advance \tempcounta by 60
\puthmorphism(\xpos,\ypos)[#2`#3`#5]{#7}{\arrowtypeb}{#9}
\advance\ypos by \tempcounta
\puthmorphism(\xpos,\ypos)[\phantom{#2}`\phantom{#3}`#4]{#7}{\arrowtypea}{#9}
\advance\ypos by -\tempcounta \advance\ypos by -\tempcounta
\puthmorphism(\xpos,\ypos)[\phantom{#2}`\phantom{#3}`#6]{#7}{\arrowtypec}{#9}
}}

\def\putarc(#1,#2)[#3`#4`#5]#6#7#8{{%
\xpos #1
\ypos #2
\width #6
\arrowlength #6
\putbox(\xpos,\ypos){#3\vphantom{#4}}%
{\advance \xpos by\arrowlength
\putbox(\xpos,\ypos){\vphantom{#3}#4}}%
\horsize{\tempcounta}{#3}%
\horsize{\tempcountb}{#4}%
\divide \tempcounta by2
\divide \tempcountb by2
\advance \tempcounta by30
\advance \tempcountb by30
\advance \xpos by\tempcounta
\advance \arrowlength by-\tempcounta
\advance \arrowlength by-\tempcountb
\halflength=\arrowlength \divide\halflength by 2
\divide\arrowlength by 5
\put(\xpos,\ypos){\bezier{\arrowlength}(0,0)(50,50)(\halflength,50)}
\ifnum #7=-1 \put(\xpos,\ypos){\vector(-3,-2)0} \fi
\advance\xpos by \halflength
\put(\xpos,\ypos){\xpos=\halflength \advance\xpos by -50
   \bezier{\arrowlength}(0,50)(\xpos,50)(\halflength,0)}
\ifnum #7=1 {\advance \xpos by
   \halflength \put(\xpos,\ypos){\vector(3,-2)0}} \fi
\advance\ypos by 50
\vertsize{\tempcounta}{#5}%
\divide\tempcounta by2
\advance \tempcounta by20
\if a#8 %
   \advance \ypos by\tempcounta
   \putbox(\xpos,\ypos){#5}%
\else
   \advance \ypos by-\tempcounta
   \putbox(\xpos,\ypos){#5}%
\fi
}}

\makeatother

\usepackage{makeidx}  

\makeatother

\usepackage{amsmath}
\usepackage{qtree}
\usepackage{dsfont}
\usepackage{stmaryrd}
\usepackage{enumerate}
\usepackage{hyperref}

\usepackage{lscape}

\newtheorem{theorem}{Theorem}[section]
\newtheorem{lemma}[theorem]{Lemma}
\newtheorem{corollary}[theorem]{Corollary}

\newtheorem{proposition}[theorem]{Proposition}

\def\true{{\bf true} }
\def\false{{\bf false} }

\def\stl{{{\bf st}^l}}
\def\str{{{\bf st}^r}}

\def\pu{{{\bf pile'up}}}
\def\pul{{{\bf pile'up}^l}}
\def\pur{{{\bf pile'up}^r}}

\def\cps{{{\bf CPS}}}
\def\cpsl{{{\bf CPS}^l}}
\def\cpsr{{{\bf CPS}^r}}

\def\trl{{{\bf TR}^l}}
\def\trr{{{\bf TR}^r}}

\def\eps{{{\bf eps}}}
\def\epsl{{{\bf eps}^l}}
\def\epsr{{{\bf eps}^r}}

\def\mos{{{\bf mos}}}
\def\mosl{{{\bf mos}^l}}
\def\mosr{{{\bf mos}^r}}
\def\strat{{{\bf strat}}}

\def\lk{\langle}
\def\rk{\rangle}

\newcommand{\lra}{\longrightarrow}
\newcommand{\hra}{\hookrightarrow}

\newcommand{\ra}{\rightarrow}
\newcommand{\la}{\leftarrow}

\newcommand{\Ra}{\Rightarrow}

\def\la{{\lambda}}

\def\bT{{\bf T}}

\def\bt{{\bf t}}

\def\cC{{\cal C}}

\def\cF{{\cal F}}

\def\cP{{\cal P}}

\def\cR{{\cal R}}

\def\cU{{\cal U}}

\def\Q{{\cal Q}}
\def\vect[#1]{\vec{#1}}

\begin{document}

\title{Scope ambiguities, monads and strengths}

\author{Justyna Grudzi\'{n}ska \and Marek Zawadowski}


\maketitle

\begin{abstract}
In this paper, we will discuss three semantically distinct scope assignment strategies: traditional movement strategy (\cite{may}, \cite{montague}, \cite{cooper}), polyadic approach (\cite{may85}, \cite{keenan87}, \cite{keenan92}, \cite{BZ}, \cite{benthem}), and continuation-based approach (\cite{barker}, \cite{groote}, \cite{barker:shan}, \cite{kiselyov:shan}, \cite{bekki}). As a generalized quantifier on a set $X$ is an element of $\cC(X)$, the value of the continuation monad $\cC$ on $X$, in all three approaches QPs are interpreted as $\cC$-computations. The main goal of this paper is to relate the three strategies to the computational machinery connected to the monad $\cC$ (strength and derived operations). As will be shown, both the polyadic approach and the continuation-based approach make heavy use of monad constructs. In the traditional movement strategy, monad constructs are not used but we still need them to explain how the three strategies are related and what can be expected of them wrt handling scopal ambiguities in simple sentences.

\end{abstract}

\section{Scope ambiguities}

Multi-quantifier sentences have been known to be ambiguous with different readings corresponding to how various quantifer phrases (QPs) are semantically related in the sentence. For example,
\begin{enumerate}[(1)]
  \item Every girl likes a boy
\end{enumerate}
admits of the subject wide scope reading ($S > O$) where each girl likes a potentially different boy, and the object wide scope reading ($O > S$) where there is one boy whom all the girls like. As the number of QPs in a sentence increases, the number of distinct readings also increases. Thus a simple sentence with three QPs admits of six possible readings, and in general a simple sentence with $n$ QPs will be (at least) $n!$ ways ambiguous (we only consider readings where QPs are linearly ordered - what we will call asymmetric readings).

In this paper, we will discuss three semantically distinct scope assignment strategies
\begin{description}
  \item[Strategy A:] Traditional movement strategy (\cite{may}, \cite{montague}, \cite{cooper}).
  \item[Strategy B:] Polyadic approach (\cite{may85}, \cite{keenan87}, \cite{keenan92}, \cite{BZ}, \cite{benthem}).
  \item[Strategy C:] Continuation-based approach (\cite{barker}, \cite{groote}, \cite{barker:shan}, \cite{kiselyov:shan}, \cite{bekki}).
\end{description}
Scope assignment strategies can be divided into two families: movement analyses (Strategies A and B) and in situ analyses (Strategy C). Strategy A has been implemented in various ways using May's QR (\cite{may}), Montague's Quantifying In Rule (\cite{montague}), Cooper's Storage (\cite{cooper}). Strategy B involving polyadic quantification has been first introduced in the works of May (\cite{may85}), Keenan (\cite{keenan87}), Zawadowski (\cite{BZ}) and van Benthem (\cite{benthem}). The most recent Strategy C involves continuations and has been first proposed in the works of Barker (\cite{barker}) and de Groote (\cite{groote}), and then further developed and modified in the works of Barker and Shan (\cite{barker:shan}), Kiselyov and Shan (\cite{kiselyov:shan}) and Bekki and Asai (\cite{bekki}). The continuation-based strategies can be divided into two groups: those that locate the source of scope-ambiguity in the rules of semantic composition and those that attribute it to the lexical entries for the quantifier words. In this paper, we only consider operation-based approaches (as in \cite{barker}). As a generalized quantifier on a set $X$ is an element of $\cC(X)$, the value of the continuation monad $\cC$ on $X$, in all three approaches QPs are interpreted as $\cC$-computations. The main goal of this paper is to relate the three scope assignment strategies to the computational machinery connected to the monad $\cC$ (strength and derived operations). As will be shown, Strategies B and C make heavy use of monad constructs. In Strategy A, monad constructs are not used but we still need them to explain how the three strategies are related and what can be expected of them wrt handling scopal ambiguities in simple sentences.

\section{Monads and strenghts}

For unexplained notions related to category theory, we refer the reader to standard textbooks on category theory. We shall be exclusively working in the cartesian closed category of sets $Set$. The category $Set$ of sets has sets as objects. A morphism in $Set$ from an object (set) X to an object (set) Y is a function $f: X\ra Y$ from $X$ to $Y$.

\subsection{Monads}

A {\em monad} on $Set$ is a triple $(T,\eta,\mu)$ where $T: Set \lra Set$ is an endofunctor (the underlying functor of the monad), $\eta : 1_{Set} \lra T$ and $\mu: T^2 \lra T$ are natural transformations (first from identity functor on $Set$ to $T$, second from the composition of $T$ with itself to $T$)
making the following diagrams
\begin{center}
\xext=2600 \yext=600 \adjust[`I;I`;I`;`I]
\begin{picture}(\xext,\yext)(\xoff,\yoff)
\settripairparms[1`1`1`1`1;500]

\putVtrianglepair(0,0)[T`T^2`T`T;\eta_T`T(\eta)`1_T`\mu`1_T ]
\putmorphism(500,500)(1,0)[\phantom{T^2}`T`T(\eta)]{500}{-1}a

\setsqparms[1`1`1`1;800`500]
\putsquare(1600,0)[T^3`T^2`T^2`T;\mu_T`T(\mu)`\mu`\mu]
\end{picture}
\end{center}
commute. $\eta$ and $\mu$ are often referred to as {\em unit} and {\em multiplication} of the monad $T$, respectively. These diagrams express the essence of the algebraic calculations. We shall explain their meaning while describing the list monad below.

Monads can serve many different purposes. Here, we think of a monad as a device to extend the notion of computation. We think of a function $f: X\lra Y$ as a computation that, when given an element of $x$, provides (computes) an element $f(x)$ of $Y$. Then the function $f:X\lra T(Y)$ can be thought of as a computation that, when given an element $x$ in $X$, provides a computation in $T(Y)$ that might, in principle, evaluate to an element of $Y$. We shall illustrate the concept on some examples below, before we focus on the continuation monad - the main notion of computation considered in this paper.

{\em Examples of monads.}
\begin{enumerate}
  \item {\em Identity monad} is the simplest possible monad but not very interesting. In this case the functor $T$ and the natural transformations $\eta$ and $\mu$ are identities.
      For this monad, the notion of a $T$-computation in $X$ is just an element of $X$, as the function $f: X\lra T(Y)$ is just $f:X\lra Y$.
  \item {\em Maybe monad} is the simplest non-trivial  monad. The functor $T$ associates to every set $X$ the set $T(X) =X+\{ \bot\}$ (the disjoint sum of $X$ and singleton $\{ \bot\}$), and to every function $f:X\lra Y$ a function $T(f): T(X) \lra T(Y)$ such that, for $x\in T(X)$,
      \[ T(f)(x)  = \left\{ \begin{array}{ll}
        x & \mbox{if $x\in X$}  \\
        \bot & \mbox{if $x=\bot$}
                                    \end{array}
                \right. \]
       So $T$ adds to $X$ an additional element $\bot$, called {\em bottom} or {\em nothing}.
       The component at $X$ of natural transformation $\eta$ is a function $\eta_X: X \lra X+\{ \bot\}$ such that $\eta_X(x)=x$, i.e. it sends $x$ to the same $x$ but in the set $X+\{ \bot\}$.  The component at $X$ of natural transformation $\mu$ is a function $\mu_X: X+\{ \bot, \bot' \} \lra X+\{ \bot\}$ such that, for $x\in  X+\{ \bot, \bot' \}$,
       \[ \mu_X(x)  = \left\{ \begin{array}{ll}
        x & \mbox{if $x\in X$}  \\
        \bot & \mbox{if $x=\bot$ or $x=\bot'$}
                                    \end{array}
                \right. \]
       i.e. it sends $x$ in $X$ to the same $x$, and two bottoms $\bot$ and $\bot'$ in $T^2(X)$ to the only bottom $\bot$ in $T(X)$.

       For this monad, the notion of a $T$-computation in $X$ consists of elements of $X$ and an additional computation $\bot$ that says that we do not get a value in $X$. The function $f: X\lra T(Y)$ is just a partial function $f:X{\hra} Y$. So this monad allows to treat partial computations as total.

  \item {\em Exception monad} is still less trivial than maybe monad. We are given a fixed set of exceptions $E$ and, for a set $X$, the monad functor is $T(X)=X+E$, i.e. the disjoint union of $X$ and $E$. If $E$ is empty, it is the identity monad; if $E$ is a singleton, then it is a maybe monad; otherwise is it like maybe monad but with many options for nothingness.
  \item {\em List monad} or {\em monoid monad} is still more interesting than the previous monad and we shall work it out in detail. It is not needed for the applications in the paper but it provides some intuitions before we move to the continuation monad. To any set $X$ the list monad functor associates the set $T(X)$ of (finite) words over $X$ (treated as an alphabet). This includes the empty word $\varepsilon$. To a function $f:X\lra Y$ the functor $T$ associates the function $T(f): T(X) \lra T(Y)$ sending the word $x_1,x_2,\ldots, x_n$ over $X$ to the word $f(x_1),f(x_2),\ldots, f(x_n)$ over $Y$. The component at $X$ of natural transformation $\eta$ is a function $\eta_X: X \lra T(X)$ such that $\eta_X(x)=x$, i.e. it sends (the letter) $x$ to the one letter word $x$ in $T(X)$.

      The component at $X$ of natural transformation $\mu$ is a function $\mu_X: T^2(X) \lra T(X)$. Note that $T^2(X)=T(T(X))$ is the set of words whose letters are words over the alphabet $X$. Thus it can be thought of as a list of lists. $\mu_X$ applied to such a list of lists flattens it to the single list. A three letter word $t=(x_1,x_2),(x_3,x_4,x_5),\varepsilon$ is a typical element of $T^2(X)$. The result of flattening $T$ is the list $\mu_X(T)=x_1,x_2,x_3,x_4,x_5$ in $T(X)$. We can think of such a word $w$ as a term/word/computation $u=y_1,y_2,y_3$ in which we intend to substitute the term $v_1=x_1,x_2$ for variable $y_1$, the term $v_2=x_3,x_4,x_5$ for variable $y_2$, and the term $v_3=\varepsilon$ for variable $y_3$, i.e. $u{[}y_1\backslash v_1,y_2\backslash v_2,y_3\backslash v_3{]}$. Now the multiplication $\mu$ can be thought of as an actual substitution. With this interpretation one can understand the intuitions behind the monad diagrams. In the left triangle, an element of $T(X)$, say $x_1,x_2,x_3$, is mapped through $\eta_{T(X)}$ to single letter word $(x_1,x_2,x_3)$ and $\mu_X$ flattens it back to $x_1,x_2,x_3$, as required for the triangle to commute. In other words, the substitution  $y[y\backslash v]$ results in $v$. In the right triangle, the map $T(\eta_x)$ sends, say $x_1,x_2,x_3$, to the letter word $(x_1),(x_2),(x_2)$ with each letter being a single letter word. Thus again flattening such a list gives $x_1,x_2,x_3$ back, as required. In other words the substitution  $y_1,y_2,y_3{[}y_1\backslash x_1,y_1\backslash x_2,y_3\backslash x_3{]}$ results in $x_1,x_2,x_3$.
      Thus the above triangles ensure that if we substitute with either term being a variable, then we get the expected result. The commutation of the square diagram, in this case, expresses the fact if we have a list of lists of lists and we flatten it in two different ways, once starting with the upper two levels of lists and the other time starting with the lower two levels of lists, and then we flatten the results again to get the ordinary lists over $X$ in $T(X)$, these lists coincide. On a more conceptual level, this square expresses the fact that the evaluation commutes with substitution. In this sense these diagrams capture the essence of all algebraic calculations.

      For this monad, the notion of a $T$-computation in $X$ consists of words over $X$ to be evaluated/computed in a monoid when elements of $X$ will be (interpreted) in a monoid. The function $f: X\lra T(Y)$ is just a function $f:X\lra T(Y)$ sending elements of $X$ to words over $Y$. So this monad allows for a list of values, for a given input.

  \item See next subsection for unexplained notation.

  {\em (Covariant) power-set monad} sends set $X$ to power-set $\cP(X)$ and a function $f:X\lra Y$ to the image function $\cP(f)=\vec{f}:\cP(X)\lra\cP(Y)$ such that, for $h:X\ra \bt$, $\vec{f}(h)$ is an image of $U$ under the function $f$, i.e. for $y\in Y$
  \[ \vec{f}(h)(y) = \bigvee_{x\in X,\; f(x)=y} h(x).  \]
  The unit $\eta_X: X \lra \cP(X)$ embeds $x$ in $X$ to the characteristic function $x$, i.e.
  \[ \eta_X(x)(x')  = \left\{ \begin{array}{ll}
        \true & \mbox{if $x=x'$}  \\
        \false & \mbox{otherwise.}
                                    \end{array}
                \right. \]
  The multiplication $\mu_X: \cP^2(X)\ra \cP(X)$ sums elements of elements of elements, i.e. for $H\in \cP^2(X)$ and $x\in X$, $\mu_X(H)(x):X\ra \bt$ is a function given by
    \[ \mu_X(H)(x)  = \bigvee_{h\in \cP(X)}\; H(h)\wedge h(x). \]

This monad also allows for a set of values, for a given input.

  \end{enumerate}

\subsection{Notation}

      Before we explain the notion of computation coming with the continuation monad, we restate the monad in a more functional way. To do this, we need to introduce some notation. As $Set$ is a cartesian closed category, it is customary to denote functions between sets using $\lambda$ notation. One can think of it as if we were to work in the internal language of $Set$, i.e. $\lambda$ theory where all functions have their names represented. For sets $X$ and $Y$, we shall use $X\times Y$ to denote the binary product of $X$ and $Y$ and $X\Ra Y$ to denote the set of functions from $X$ to $Y$. As it is customary, we associate $\Ra$ to the right, i.e. $X\Ra Y\Ra Z$ means $X\Ra (Y\Ra Z)$ and this set is naturally bijective with $(X\times Y)\Ra Z$. If we have a function $$f: X\times Y \lra Z,$$ then by
      $$ \lambda y_{:Y}. f : X \lra Y\Ra Z $$
      we denote its exponential adjunction, i.e. the function from $X$ to the set of functions $Y\Ra Z$ such that, for an element $x\in X$, $\lambda y_{:Y}. f(x)$  is a function from $Y$ to $Z$ such that, for an element $y\in Y$, $(\lambda y_{:Y}. f)(x)(y)$ is by definition equal $f(x,y)$. Note that in the expression $(\lambda y_{:Y}. f)(x)(y)$ the first occurrence of $y$ is an occurrence of a variable (as it is part of the name of a function), whereas the second occurrence of $y$ in this expression denotes an element of the set $Y$.

$\pi_i$ will denote the projection on $i$-component from the product. Any function $\sigma: \{1,\ldots m\} \ra \{ 1,\ldots, n\}$ induces a generalized projection denoted
\[ \pi_\sigma =\lk \pi_{\sigma(1)},\ldots, \pi_{\sigma(m)} \rk: X_1\times \ldots\times X_n\lra X_1\times \ldots\times X_m. \]
 We will use this notation mainly for $\sigma$'s being bijections, i.e. when $\pi_\sigma$ is just permutation of the component for the product.

We have a distinguished set of truth values $\bt=\{ \true, \false\}$. We shall use the usual (possibly infinitary) operations on this set. For a set $X$, we put $\cP(X)=X\Ra \bt$, i.e. the (functional) powerset of $X$.

\subsection{Continuation monad}
  {\em Continuation monad}, the most important for us, denoted $\cC$, has some similarities to the power-set monad but it also differs in a substantial way. At the level of objects, it is just twice iterated power-set construction, i.e. for set $X$, $\cC(X)=\cP^2(X)$, but at the level of morphisms, it is an inverse image of an inverse image, i.e., function $f:X\ra Y$ induces an inverse image function between powersets \[ \cP(f)=f^{-1}: \cP(Y)\ra \cP(X) \]
       \[ h\mapsto h\circ f, \hskip 2mm \cP(f)= \la h_{:\cP(Y)}.\la x_{:X}. h(f\, x) \]
  Taking again an inverse image function, we have
      \[ \cC(f)=\cP(f^{-1}) : \cC(X) \ra \cC(Y) \]
  \[ Q \mapsto Q\circ f^{-1} , \hskip 2mm \cC(f)(Q)= \lambda h_{:\cP(Y)}.  Q(\lambda x_{:X}. h( f\, x))   \]
  for $Q\in \cC(X)$.\\
  The unit $\eta_X : X\ra \cC(X)$ is given by
  \[ \eta_X(x)= \lambda h_{:\cP(X)}. h(x).\]
  for $x\in X$.\\
  The multiplication $\mu_X: \cC^{2}(X)\lra \cC(X)$ can be explained in terms of $\eta$
  \[ \mu_X= (\eta_{\cP(X)})^{-1}: \cP^{4}(X)\lra \cP^{2}(X).\]
  In other words, $\mu_X(\cF):\cP(X)\ra \bt$ is a function such that
  \[  \mu_X(\cF)(h) = \cF(\eta_{\cP(X)}(h))\]
   for $\cF :\cP^{3}(X)\ra \bt$ and $h: X\ra \bt$.\\
  In $\lambda$-notation, we write
  \[\mu_X(\cF)(h)=\cF(\la D_{: \cC(X)}.D(h)).\]

  Now we can look at the notion of computation related the continuation monad. Consider the function
  \[ f : X\lra \cC(Y) .\]
  By exponential adjunction (uncurrying) it corresponds to a function
  \[ f' : \cP(Y)\times X\lra \bt \]
  and again by exponential adjunction (currying) it corresponds to a function
  \[ f'' : \cP(Y)\lra \cP(X). \]
Thus a $\cC$-computation from $X$ to $Y$ is a function that sends functions from $\cP(Y)=Y\Ra\bt$ to functions in  $\cP(X)$. So instead of having a direct answer for a given element $x\in X$ what is the value $f(x)$ in $Y$, we are given for every continuation function $c: Y\lra \bt$ a value in the answer type $\bt$ that could be thought of as $c(f(x))$ (if there were an element in $Y$ that could be reasonably called $f(x)$). We can draw the picture illustrating the situation
\begin{center} \xext=1200 \yext=300
\begin{picture}(\xext,\yext)(\xoff,\yoff)
  \putmorphism(0,50)(1,0)[X`Y`f?]{600}{1}r
  \putmorphism(600,50)(1,0)[\phantom{Y}`\bt`c]{600}{1}r
    \put(630,305){$f(c)$}
     \put(0,150){\line(0,1){100}}
  \put(0,250){\line(1,0){1200}}

   \put(1200,250){\vector(0,-1){100}}
 \end{picture}
\end{center}
Instead of `procedure' $f?$ computing $y$'s from $x$'s (that we don't have), we provide a continuation $f(c)$ for any continuation (of the computation) $c$. If $f?$ would be indeed a genuine function $f?: X\lra Y$, then $f(c)$ would be the composition $c\circ f?$.

\subsection{Bi-strong monads}

As the notion of strength is new in this context, we shall briefly recall its history. There are three manifestations of strength on a functor. Historically, the first one was the notion of enrichment of a functor (c.f. \cite{eilenberg}). The tensorial strength (i.e., natural transformation of a kind $X\otimes T(Y)\lra T(X\otimes Y)$  used in this paper) was introduced in  \cite{kock1} and further developed in  \cite{kock4}. The cotensorial strength (i.e., natural transformation of a kind $T(X\Ra Y)\ra X\Ra T(Y)$) introduced in \cite{kock3} also proved useful in some contexts. In symmetric monoidal closed categories these concepts are equivalent, (c.f. \cite{kock3}).

As it was noticed in \cite{moggi}, a monad, in order to have a well behaved notion of computation, has to be strong. Fortunately, all monads on $Set$ are strong. More precisely, all monads on set can be canonically equipped with two strengths, left and right, and moreover these strengths are compatible in a precise technical sense. This additional structure on the continuation monad will be essential when we shall analyze the meaning of multiple quantified sentences.

Let $(T,\eta,\mu)$ be a monad on $Set$. The {\em left strength} is a natural transformation with components
\[  \stl_{X,Y} : T(X)\times Y \lra T(X\times Y) \]
for sets $X$ and $Y$, making the diagrams
\begin{center}
\xext=1100 \yext=600 \adjust[`I;I`;I`;`I]
\begin{picture}(\xext,\yext)(\xoff,\yoff)
\settriparms[1`1`-1;550]
  \putVtriangle(0,0)[T(X)\times Y\times Z`T(X\times Y\times Z)`T(X\times Y)\times Z;\stl_{X,Y\times Z}`\stl_{X,Y}\times 1` \stl_{X\times Y, Z}]
\end{picture}
\end{center}
and
\begin{center}
\xext=1500 \yext=1100
\begin{picture}(\xext,\yext)(\xoff,\yoff)
\putmorphism(0,1000)(0,-1)[X\times Y`\phantom{T(X)\times Y}`\eta_X \times 1]{500}{1}l

\putmorphism(50,1000)(2,-1)[\phantom{X\times Y}`\phantom{T(X\times Y)}`\eta_{X\times Y}]{1000}{1}r

 \setsqparms[1`-1`0`1;1200`500]
\putsquare(0,0)[T(X)\times Y`T(X\times Y)`T^2(X)\times Y`T(T(X)\times Y);{\stl}_{X,Y}`\mu_X\times 1``\stl_{T(X),Y}]

\putmorphism(1250,500)(2,-1)[\phantom{T(X\times Y)}`\phantom{T^2(X\times Y)}`\mu_{X\times Y}]{1000}{-1}r
\putmorphism(1200,0)(1,0)[\phantom{T(X\times T(Y))}`T^2(X\times Y)`T({\stl}_{X\times Y})]{1200}{1}b

\end{picture}
\end{center}
commute.

The {\em right strength} is a natural transformation with components
\[  \str_{X,Y} : X\times T(Y) \lra T(X\times Y) \]
for sets $X$ and $Y$, making the diagrams
\begin{center}
\xext=1100 \yext=600 \adjust[`I;I`;I`;`I]
\begin{picture}(\xext,\yext)(\xoff,\yoff)
\settriparms[1`1`-1;550]
  \putVtriangle(0,0)[X\times Y\times T(Z)`T(X\times Y\times Z)`X\times T(Y\times Z);\str_{X\times Y,Z}`1\times \str_{Y,Z}` \str_{X,Y\times Z}]
\end{picture}
\end{center}
and
\begin{center}
\xext=1500 \yext=1100
\begin{picture}(\xext,\yext)(\xoff,\yoff)
\putmorphism(0,1000)(0,-1)[X\times Y`\phantom{X\times T(Y)}`1\times \eta_Y]{500}{1}l

\putmorphism(50,1000)(2,-1)[\phantom{X\times Y}`\phantom{T(X\times Y)}`\eta_{X\times Y}]{1000}{1}r

 \setsqparms[1`-1`0`1;1200`500]
\putsquare(0,0)[X\times T(Y)`T(X\times Y)`X\times T^2(Y))`T(X\times T(Y));{\str}_{X,Y}`1\times \mu_Y``\str_{X,T(Y)}]

\putmorphism(1250,500)(2,-1)[\phantom{T(X\times Y)}`\phantom{T^2(X\times Y)}`\mu_{X\times Y}]{1000}{-1}r
\putmorphism(1200,0)(1,0)[\phantom{T(X\times T(Y))}`T^2(X\times Y)`T({\str}_{X\times Y})]{1200}{1}b

\end{picture}
\end{center}
commute.

The monad $(T,\eta,\mu)$ on $Set$ together with two natural transformations $\stl$ and $\str$ of right and left strength is a {\em bi-strong monad} if, for any sets $X$, $Y$, $Z$, the square
\begin{center}
\xext=1200 \yext=600 \adjust[`I;I`;I`;`I]
\begin{picture}(\xext,\yext)(\xoff,\yoff)
\setsqparms[1`1`1`1;1200`500]
\putsquare(00,0)[X\times T(Y)\times Z`X\times T((Y\times Z)`T(X\times Y)\times Z`T(X\times Y \times Z);1_X\times \stl_{Y,Z}`\str_{X,Y}\times 1_Z`\str_{X,Y\times Z}`\stl_{X\times Y, Z}]
\end{picture}
\end{center}
commutes.

As we already mentioned, each monad $(T,\eta,\mu)$ on $Set$ is bi-strong. We shall define the right and left strength. Fix sets $X$ and $Y$.
For $x\in X$ and $y\in Y$, we have functions
\[ l_{y} : X\lra X\times Y, \hskip 5mm{\rm and} \hskip 5mm  r_{x} : Y\lra X\times Y,\]
such that
\[ l_{y}(x) = \lk x,y\rk,  \hskip 5mm{\rm and} \hskip 5mm r_{x}(y) = \lk x,y\rk. \]
The left and right strength
\[ \stl_{X,Y}: T(X)\times Y \lra T(X\times Y) \hskip 5mm{\rm and} \hskip 5mm \str_{X,Y}: X\times T(Y) \lra T(X\times Y) \]
are given for $x\in X$, $s\in T(X)$, $y\in Y$ and $t\in T(Y)$ by
\[ \stl_{X,Y}(s,y)= T(l_y)(s)   \hskip 5mm{\rm and} \hskip 5mm  \str_{X,Y}(x,t)= T(r_x)(t), \]
respectively. We drop indices $_{X,Y}$ when it does not lead to confusion.

It is not difficult to verify that the above defines left ($\stl$) and right ($\str$) strength on the monad $T$ and since, for any $x\in X$ and $z\in Z$, the square
\begin{center}
\xext=1200 \yext=600 \adjust[`I;I`;I`;`I]
\begin{picture}(\xext,\yext)(\xoff,\yoff)
\setsqparms[1`1`1`1;1200`500]
\putsquare(00,0)[Y`X\times Y`Y\times Z`X\times Y \times Z;r_x`l_z`l_z`r_x]
\end{picture}
\end{center}
commutes, they are compatible and make the monad $T$ bi-strong. Note that these strengths are related by the following diagram
\begin{center}
\xext=1200 \yext=600 \adjust[`I;I`;I`;`I]
\begin{picture}(\xext,\yext)(\xoff,\yoff)
\setsqparms[1`1`-1`1;1200`500]
\putsquare(00,0)[T(X)\times Y `T(X\times Y)`Y\times T(X)`T(Y\times X);\stl_{X,Y}`T(\lk \pi_2,\pi_1\rk)`\lk \pi_2,\pi_1\rk`\str_{Y,X}]
\end{picture}
\end{center}

{\em Examples of strength on monads in $Set$.}
\begin{enumerate}
  \item Maybe monad. The left strength $\stl_{X,Y} : (X+\{\bot\}) \times Y \lra (X \times Y) +\{\bot\}$ is given by
  \[ \stl(x,y)  = \left\{ \begin{array}{ll}
        \bot & \mbox{if $x=\bot$}  \\
        \lk x,y\rk & \mbox{otherwise.}
                                    \end{array}
                \right. \]
      Right strength is similar.
  \item List monad. The left strength $\stl : T(X) \times Y \lra T(X \times Y)$ is given by
  \[ \stl(\vec{x},y)  = \left\{ \begin{array}{ll}
        \varepsilon & \mbox{if $\vec{x}=\varepsilon$}  \\
        \lk x_1,y\rk, \ldots,\lk x_n,y\rk  & \mbox{if $\vec{x}=x_1,\ldots, x_n$.}
                                    \end{array}
                \right. \]
      Right strength is similar.
  \item Continuation monad. We shall describe the strength morphisms by lambda terms. The left strength is
  $$\stl=\la N_{:\cC(X)}.\la y_{:Y}.  \la c_{:\cP(X\times Y)}. N(\la x_{:X}. c(x,y)): \cC(X)\times Y \lra \cC(X\times Y)$$
  and the right strength is
  $$\str=\la x_{:X}. \la M_{:\cC(Y)}. \la c_{: \cP(X\times Y)}. M(\lambda y_{:Y}.c(x,y)): X\times \cC(Y) \lra \cC(X\times Y).$$

\end{enumerate}

\subsection{Combining computations in arbitrary monad $T$ on $Set$}

Using both strengths, we can define two {\em pile up} natural transformations, the left  and right. For any sets $X$ and $Y$, the {\em left pile up} $\pul_{X,Y}$ is defined from the diagram
\begin{center}
\xext=1200 \yext=600 \adjust[`I;I`;I`;`I]
\begin{picture}(\xext,\yext)(\xoff,\yoff)
\setsqparms[1`1`-1`1;1200`500]
\putsquare(00,0)[T(X)\times T(Y)`T(X\times Y)`T(X\times T(Y))`T^2(X\times Y);{\pul}_{X,Y}`\stl_{X,T(Y)}`\mu_{X\times Y}`T(\str_{X,Y})]
\end{picture}
\end{center}
In the above diagram, the function ${\pul}_{X,Y}$ is defined as a composition of three operations: the first is taking the $T$-computation on $X$ `outside' to be a computation on $X\times T(Y)$, the second is taking the $T$-computation on $Y$ `outside' to be a $T$-computation on $X\times Y$. In this way, we have  $T$-computations coming from $X$ on $T$-computations coming from $Y$ on $X\times Y$. Now the last morphism $\mu_{X\times Y}$ flattens these two levels to one, i.e. the $T$-computation on $T$-computations to $T$-computations.

The {\em right pile up} ${\pur}_{X,Y}$ is defined from the diagram
\begin{center}
\xext=1200 \yext=600 \adjust[`I;I`;I`;`I]
\begin{picture}(\xext,\yext)(\xoff,\yoff)
\setsqparms[1`1`-1`1;1200`500]
\putsquare(00,0)[T(X)\times T(Y)`T(X\times Y)`T(T(X)\times Y)`T^2(X\times Y);\pur_{X,Y}`\str_{T(X),Y}`\mu_{X\times Y}`T(\stl_{X,Y})]
\end{picture}
\end{center}
This operation takes out the $T$-computations in a reverse order and so they pile up in the other way.

If these {\em pile up} operations agree for all sets $X$ and $Y$, the monad is called {\em commutative}. On our list of monads, identity, maybe and covariant power-set monads are commutative. The exception, list and continuation monads are not commutative. Most monads, including the continuation monad $\cC$, are not commutative. It should be noticed that even if the monad $T$ is not commutative, both lift morphisms agree on pairs in which at least one component comes from the actual value (not an arbitrary $T$-computation). In other words, the functions
\begin{center} \xext=1400 \yext=400
\begin{picture}(\xext,\yext)(\xoff,\yoff)
  \putmorphism(0,200)(1,0)[\phantom{T(X_1)\times  T(X_2)}`\phantom{T(X_1\times X_2)}`{\pul}_{X_1,X_2}]{1400}{1}a
  \putmorphism(0,150)(1,0)[T(X_1)\times  T(X_2)`T(X_1\times X_2)`]{1400}{0}a
  \putmorphism(0,100)(1,0)[\phantom{T(X_1)\times  T(X_2)}`\phantom{T(X_1\times X_2)}`{\pur}_{X_1,X_2}]{1400}{1}b
 \end{picture}
\end{center}
are equalized by both
\begin{center} \xext=1200 \yext=100
\begin{picture}(\xext,\yext)(\xoff,\yoff)
  \putmorphism(0,00)(1,0)[X_1\times  T(X_2)`T(X_1)\times T(X_2)`\eta_{X_1}\times 1]{1200}{1}a
 \end{picture}
\end{center}
and
\begin{center} \xext=1200 \yext=100
\begin{picture}(\xext,\yext)(\xoff,\yoff)
  \putmorphism(0,00)(1,0)[T(X_1)\times  X_2`T(X_1)\times T(X_2)`1\times \eta_{X_2}]{1200}{1}a
 \end{picture}
\end{center}
morphisms. Both $\pul$ and $\pur$ are associative. All this is shown in the Appendix.

{\em Examples of $\pu$-operations.}

\begin{enumerate}
  \item Maybe monad. The left and right pile'up's coincide in this case, as in any commutative monad. We have
   $\pul_{X,Y}=\pur_{X,Y} : (X+\{\bot\}) \times (Y+\{\bot'\}) \lra (X \times Y) +\{\bot\}$ is given by
  \[ \stl(x,y)  = \left\{ \begin{array}{ll}
        \bot & \mbox{if $\{x,y\} \cap \{\bot,\bot'\}\neq\emptyset$}  \\
        \lk x,y\rk & \mbox{otherwise.}
                                    \end{array}
                \right. \]

  \item List monad. The left pile'up $\pul : T(X) \times T(Y) \lra T(X \times Y)$ is given by
  \[ \pul(\lk x_1,\ldots, x_n\rk, \lk y_1,\ldots, y_m\rk)  = \]
  \[ = \lk \lk x_1,y_1\rk,\lk x_1,y_2\rk,\ldots,\lk x_1,y_m\rk,\lk x_2,y_1\rk, \ldots,  \lk x_n,y_{m-1}\rk,\lk x_n,y_m\rk \]
   and the right pile'up $\pur : T(X) \times T(Y) \lra T(X \times Y)$ is given by
  \[ \pur(\lk x_1,\ldots, x_n\rk, \lk y_1,\ldots, y_m\rk)  = \]
  \[ = \lk \lk x_1,y_1\rk,\lk x_2,y_1\rk,\ldots,\lk x_n,y_1\rk,\lk x_1,y_2\rk, \ldots,  \lk x_{n-1},y_m\rk,\lk x_n,y_m\rk.\]

  \item {\em (Covariant) power-set monad}. The left and right pile'up's coincide in this case. We have
   $\pul_{X,Y}=\pur_{X,Y} : \cP(X) \times \cP(Y) \lra \cP(X \times Y)$ given by
  \[ \stl(U,V)  = U\times V\]
 for $U\in \cP(X)$ and $V\in \cP(Y)$.\\
  \item Continuation monad. Both pile'up operations
  $$\pul,\pur: \cC(X)\times \cC(Y) \lra \cC(X\times Y)$$
  can be defined, for $M\in \cC(X)$ and $N\in \cC(Y)$, by lambda terms as
     \[ \pul(M,N)=\la c_{:\cP(X\times Y)}. M(\la x_{:X}. N(\la y_{:Y} c(x,y)) \]
  and
  \[ \pur(M,N)= \la c_{:\cP(X\times Y)}. N(\la y_{:Y}. M(\la x_{:X} c(x,y)). \]
  The calculations for these operations are in the Appendix.

  Thus in the case of the continuation monad `piling up' computations one on top of the other is nothing but putting (interpretations of) quantifiers (= computations in the continuation monad) in order, either first before the second or the second before the first.
\end{enumerate}

\subsection{$T$-transforms on arbitrary monad $T$ on $Set$}

There are two (binary) $T$-transformations, right and and left. For a function $f: X\times Y \lra T(Z)$, the left $T$-transform is defined as the composition
\begin{center}
\xext=1200 \yext=600 \adjust[`I;I`;I`;`I]
\begin{picture}(\xext,\yext)(\xoff,\yoff)
\setsqparms[1`1`-1`1;1200`500]
\putsquare(00,0)[T(X)\times T(Y)`T(Z)`T(X\times Y)`T^2(Z);{{\trl}^{,T}}_{X,Y}(f)`{\pul}`\mu_{Z}`T(f)]
\end{picture}
\end{center}
and the right $T$-transform is defined as the composition
\begin{center}
\xext=1200 \yext=600 \adjust[`I;I`;I`;`I]
\begin{picture}(\xext,\yext)(\xoff,\yoff)
\setsqparms[1`1`-1`1;1200`500]
\putsquare(00,0)[T(X)\times T(Y)`T(Z)`T(X\times Y)`T^2(Z);{{\trr}^{,T}}_{X,Y}(f)`{\pur}`\mu_{Z}`T(f)]
\end{picture}
\end{center}
The most popular $\cps$-transforms are for the evaluation morphism $ ev : X \times (X \Ra Y) \rightarrow Y$ but there are also other morphisms having useful transforms.

 {\em Examples of $T$-transforms and in particular $\cps$-transforms.}

\begin{enumerate}
  \item  The evaluation map $ev : X \times (X \Ra Y) \rightarrow Y$ gives rise to application transforms
      \[{\trl}^{,T}(ev), {\trr}^{,T}(ev) :  T(X) \times T(X \Ra Y) \rightarrow T(Y).\]
      In case $T$ is the continuation monad $\cC$, they are the usual $\cps$-transforms
      $\cpsl(ev), \cpsr(ev) :  \cC(X) \times \cC(X\Ra Y) \rightarrow \cC(Y)$
      given  by
      \[ \cpsl(ev)(M,N)=\la h_{:\cP(Y)}. M(\la x_{:X}. N(\la g_{:X\Ra Y}. h(g\, x)))\]
      for $M\in \cC(X)$ and $N\in \cC(X\Ra Y)$.\\
      Right version is similar.

  \item Various epsilon maps are typically defined as maps from a product. Thus they give rise to various $T$-transforms. We list some of them below mainly to introduce notation that will be used later. The definitions are given by lambda terms.
      \begin{enumerate}
        \item Left evaluation
         $$\epsl_X= \la h_{:\cP(X)}.\la x_{:X}. h(x)  : \cP(X) \times X \rightarrow \bt;$$
        \item Right evaluation
         $$\epsr_X= \la x_{:X}.\la h_{:\cP(X)}. h(x)  : X \times \cP(X) \rightarrow \bt;$$

         \item Left partial  evaluation
        $$\epsl^{,X}_Y=\epsl_Y= \la c_{:\cP(X\times Y)}.\la y_{:Y}. \la x_{:X}. c(x,y)  : \cP(X\times Y)\times Y \rightarrow \cP(X);$$
        \item Right partial  evaluation
        $$\epsr^{,X}_Y=\epsr_Y=\la y_{:Y}. \la c_{:\cP(X\times Y)}. \la x_{:X}. c(x,y)  : Y\times \cP(X\times Y) \rightarrow \cP(X);$$
   \end{enumerate}
   \item What we call Mostowski maps are maps similar to $\eps$'es that are the algebraic counterpart of the interpretation of generalized quantifiers of Mostowski. Again, we give a definition for total and partial case.
   \begin{enumerate}
        \item  Left Mostowski
        $$\mosl_X =\la Q_{:\cC(X)}.\la c_{:\cP(X)}. Q(c): \cC(X)\times \cP(X) \rightarrow \bt;$$
  \item  Right Mostowski
         $$\mosr_X =\la c_{:\cP(X)}. \la Q_{:\cC(X)}. Q(c) :\cP(X) \times \cC(X) \rightarrow \bt;$$
        \item  Left partial Mostowski
         $$\mosl_Y =\la Q_{:\cC(Y)}.\la c_{:\cP(X\times Y)}. \la x_{:X}. Q( \la y_{:Y}. c(x,y)) : \cC(Y)\times \cP(X\times Y) \rightarrow \cP(X);$$
  \item  Right partial Mostowski
         $$\mosr_Y =\la c_{:\cP(X\times Y)}. \la Q_{:\cC(Y)}. \la x_{:X}. Q( \la y_{:Y}. c(x,y)) :\cP(X\times Y) \times \cC(Y) \rightarrow \cP(X).$$
      \end{enumerate}
\end{enumerate}

\section{Scope assignment strategies}

Using the notions connected to the continuation monad introduced above, we shall now precisely state and compare three strategies (A, B, C) for determining the meaning of multi-quantifier sentences.

\subsection{General remarks}

In each strategy, the starting point is the Surface Structure Tree of a sentence. This tree is rewritten so as to obtain Formal Structure Trees that correspond to all and only the available meanings of the sentence. Finally, we relabel those trees to obtain Computation Trees\footnote{We think of Computation Trees by analogy with mathematical  expressions, e.g.
\[ ((2-7)-8)+((12+5):7) \]
that can be represented as
\par
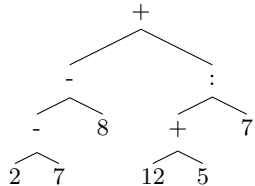
\Tree [.+ [.- [.- 2 7 ] 8 ][.: [.+ 12 5 ] 7 ]]

\medskip
\noindent i.e. a labeled binary tree where the leaves of this tree are labeled with values and the internal nodes are labeled with operations that will be applied in the computation to the values obtained from the computations of the left and right subtrees.} that provide the semantics for the sentence in each of its reading.

\vskip 2mm

\begin{center} \xext=3150 \yext=500
\begin{picture}(\xext,\yext)(\xoff,\yoff)
    \put(40,290){Surface}
    \put(20,180){Structure}
 \put(100,70){Tree}
 \put(0,50){\framebox(450,360){ }}

 \put(650,400){rewriting}
 \put(480,290){(disambiguation)}
 \put(450,250){\vector(1,0){800}}

     \put(1290,290){Formal}
    \put(1270,180){Structure}
 \put(1350,70){Tree}
 \put(1250,50){\framebox(450,360){ }}

 \put(1850,400){relabelling}
 \put(1750,290){(interpretation)}
 \put(1700,250){\vector(1,0){750}}

      \put(2490,290){Computation}
    \put(2500,180){(Semantic)}
 \put(2620,70){Tree}
 \put(2450,50){\framebox(650,360){ }}

 \end{picture}
\end{center}
\textbf{Rewriting}. Scope assignment strategies can be divided into two families: movement analyses (rewriting rules include QR, Predicate Collapsing and possibly Rotation) and in situ analyses (no rewriting rules). Below we define three rewrite rules on trees: QR Rule, Predicate Collapsing and Rotation.
\begin{itemize}
  \item QR (Quantier Raising) Rule
  \begin{itemize}
    \item applies when we have a chosen QP in a leaf of a tree;
    \item adjoins QP to S;
    \item indexes S with the variable bound by the raised QP.
    \medskip

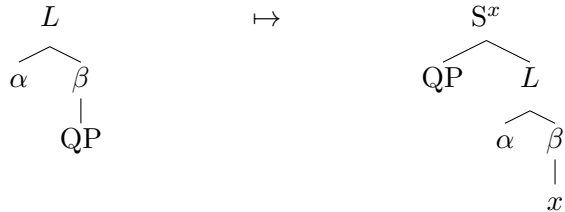
\Tree [.$L$ $\alpha$ [.$\beta$ QP ]]
\hskip 2cm $\mapsto$
\Tree [.S$^x$ QP [.$L$ $\alpha$ [.$\beta$ $x$ ]]]

\par
\medskip

($L$ - label, $\alpha, \beta$ - subtrees.)

\end{itemize}

\item Predicate Collapsing
\begin{itemize}
  \item applies when all the leaves under the node labeled S are labeled with variables (not QP's);
  \item collapses the whole subtree with the root S to a single leaf labeled with the variables $x_1, x_2, x_3$ from the leaves under S-node.
\end{itemize}

\medskip
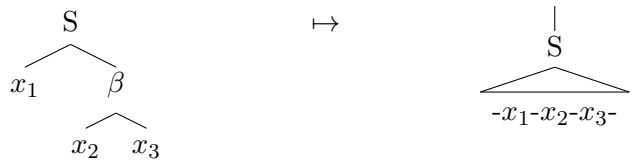
\Tree [.S [.$x$_1 ][.$\beta$ [.$x_2$ ][.$x_3$ ]]]
\hskip 2cm $\mapsto$
\Tree [\qroof{-$x_1$-$x_2$-$x_3$-}.S ]

\par
\medskip
  \item Rotation
  \begin{itemize}
  \item applies to a tree with two distinguished nodes labeled with $S$'es superscriped with some variables: the mother labeled S$^{\vec{x}}$ and its right daughter labeled S$^{\vec{y}}$;
  \item it rotates left the subtree with root labeled S$^{\vec{x}}$;
  \item the root of this subtree is labeled S$^{\vec{x}\vec{y}}$ and the (new) leftt daughter is labeled $Polyadic$.

\medskip
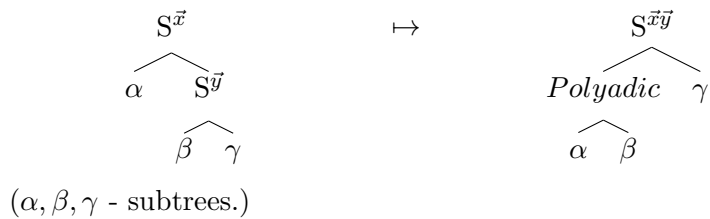
\Tree [.S$^{\vec x}$ $\alpha$ [.S$^{\vec y}$ [.$\beta$ ] [.$\gamma$ ]]]
\hskip 2cm $\mapsto$
\Tree [.S$^{\vec x \vec y}$ [.$Polyadic$ [.$\alpha$ ][.$\beta$ ]][.$\gamma$ ]]

\par
\medskip

($\alpha, \beta, \gamma$ - subtrees.)

\end{itemize}

\end{itemize}
\textbf{Relabelling}. In each scope assignment strategy, the leaves in the Computation Tree have the same labels: QPs are interpreted as $\cC$-computations, and predicates are interpreted as usual or lifted. The main difference among the three approaches consists in the shape of the Formal Structure Trees and the operations ($\eps$'es, $\mos$'es, $\pu$'es, $\cps$'es) used as labels of the inner nodes of the Computation Trees.

\subsection{Strategy A}

In the traditional movement strategy (as variously implemented in \cite{may}, \cite{montague})
\begin{itemize}
  \item Surface Structure Tree gets rewritten (disambiguated) as Formal Structure Trees (Logical Forms) via
  \begin{itemize}
  \item QR Rule;
  \item Predicate Collapsing.
\end{itemize}
  \item Formal Structure Trees (LFs) are relabelled as Computation Trees as follows
  \begin{itemize}
  \item S$^x$ (root of a subtree representing a formula) is interpreted as a suitably typed $\mos$-operation (the only operation allowed);
  \item S (leaf of a tree) is interpreted as a predicate;
  \item QP (leaf of a tree) is interpreted as a generalized quantifier $\|Q\|$ quantifying over a set $X$ (i.e. as a $\cC$-computation on $X$).
\end{itemize}
\end{itemize}
We will illustrate each strategy on the examples involving one, two and three QPs.

\vskip 2mm

\noindent \textbf{Sentence with one QP}, e.g. \textit{Every kid (most kids) entered.}

\begin{flushleft}
(A1) Surface Structure Tree

\medskip
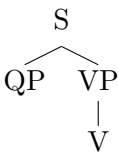
\Tree [.S [.QP ][.VP V ]]
\end{flushleft}

\begin{flushleft}
(A1)  Formal Structure Tree (LF) and the corresponding Computation Tree

\medskip
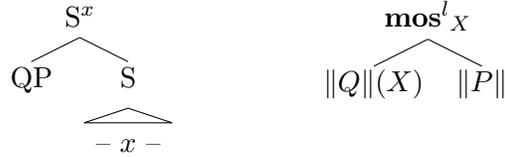
\Tree [.S$^x$ QP \qroof{-- $x$ --}.S ]
\Tree [.$\mosl_X$ $\|Q\|(X)$ $\|P\|$  ]
\par
\medskip
\end{flushleft}
The Computation Tree in (A1) gives rise to the following general map

\begin{center} \xext=600 \yext=600
\begin{picture}(\xext,\yext)(\xoff,\yoff)
  \put(-500,300){$\strat^1_A:$}
  \putmorphism(300,500)(0,-1)[\cC(X)\times\cP(X)`2`\mosl_{X}]{500}{1}r
 \end{picture}
\end{center}
In this case, there is one such map - thus Strategy A yields one reading for a sentence with one QP.

\vskip 2mm

\noindent \textbf{Sentence with two QPs}, e.g. \textit{Every girl likes a boy.}

\begin{flushleft}
(A2) Surface Structure Tree

\medskip
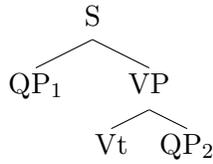
\Tree [.S [.QP_1 ][.VP [.Vt ][.QP_2 ]]]
\end{flushleft}

\begin{flushleft}
(A2) Formal Structure Tree (LF) and the corresponding Computation Tree

\medskip
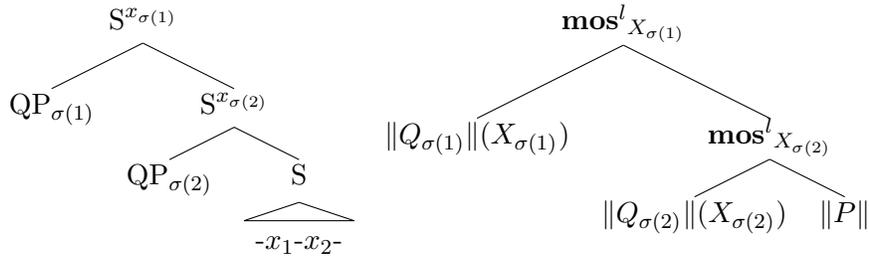
\Tree [.S$^{x_{\sigma(1)}}$ QP_{\sigma(1)} [.S$^{x_{\sigma(2)}}$ [.QP_{\sigma(2)} ] \qroof{-$x_1$-$x_2$-}.S ]]
\Tree [.$\mosl_{X_{\sigma(1)}}$  [.$\|Q_{\sigma(1)}\|(X_{\sigma(1)})$ ][.$\mosl_{X_{\sigma(2)}}$  [.$\|Q_{\sigma(2)}\|(X_{\sigma(2)})$  ] [.$\|P\|$  ]]]
\par
\medskip
\end{flushleft}
The Computation Tree in (A2) gives rise to the following general map, with $\sigma\in S_2$ (where $S_2$ is the set of permutations of the set $\{1, 2\}$)

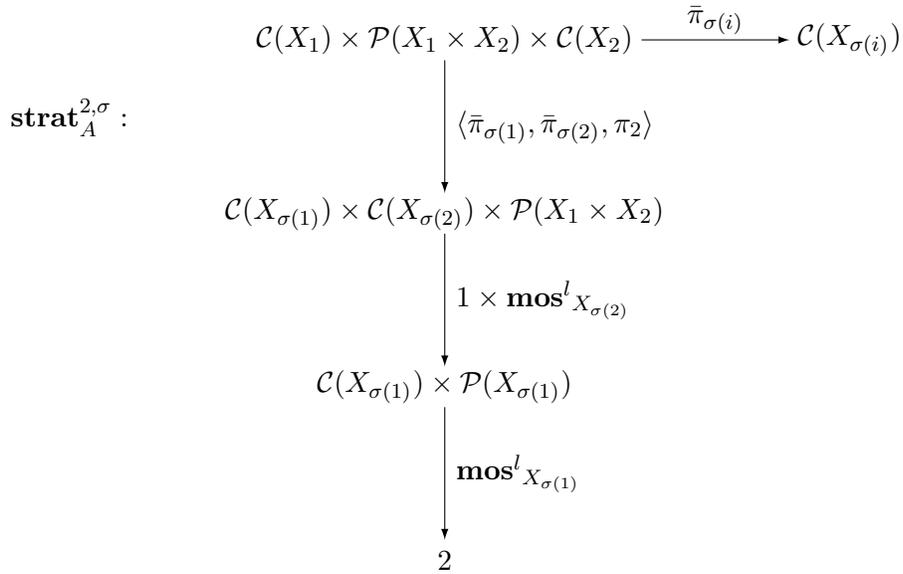
\begin{center} \xext=1000 \yext=1900
\begin{picture}(\xext,\yext)(\xoff,\yoff)
  \put(-1000,1500){$\strat^{2,\sigma}_A:$}
  \putmorphism(500,1800)(0,-1)[\cC(X_1)\times\cP(X_1\times X_2)\times\cC(X_2)`\phantom{\cC(X_1)\times\cC(X_2)\times\cP(X_1\times X_2)}`\lk \bar{\pi}_{\sigma(1)},\bar{\pi}_{\sigma(2)},\pi_2\rk]{600}{1}r

  \putmorphism(500,1800)(1,0)[\phantom{\cC(X_1)\times\cP(X_1\times X_2)\times\cC(X_2)}`\cC(X_{\sigma(i)})`\bar{\pi}_{\sigma(i)} ]{1400}{1}a

  \putmorphism(500,1200)(0,-1)[\cC(X_{\sigma(1)})\times\cC(X_{\sigma(2)})\times\cP(X_1\times X_2)`\phantom{\cC(X)\times\cC(\cP(X))}`1\times \mosl_{X_{\sigma(2)}}]{600}{1}r
  \putmorphism(500,600)(0,-1)[\cC(X_{\sigma(1)})\times\cP(X_{\sigma(1)})`2`\mosl_{X_{\sigma(1)}}]{600}{1}r
 \end{picture}
\end{center}
where $\bar{\pi}_{\sigma(i)}$ is the projection on the 1st factor if $\sigma(i)=1$ and on the 3rd factor if $\sigma(i)=2$, i.e. as it should be. This convention will be used in all similar diagrams without any further explanations.\\
There are two such maps corresponding to the two permutations $\sigma$ of $\{ 1,2\}$. These maps are different in general.
Thus Strategy A yields two (both) asymmetric readings for a sentence with two QPs.

\vskip 2mm

\noindent \textbf{Sentence with three QPs}, e.g. \textit{Some teacher gave every student most books.}

\begin{flushleft}
(A3) Surface Structure Tree

\medskip
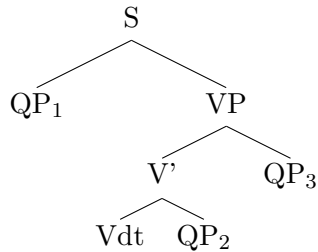
\Tree [.S [.QP_1 ][.VP [.V' [.Vdt ][.QP_2 ]][.QP_3 ]]]
\end{flushleft}

\bigskip

\begin{flushleft}
(A3) Formal Structure Tree (LF)

\medskip
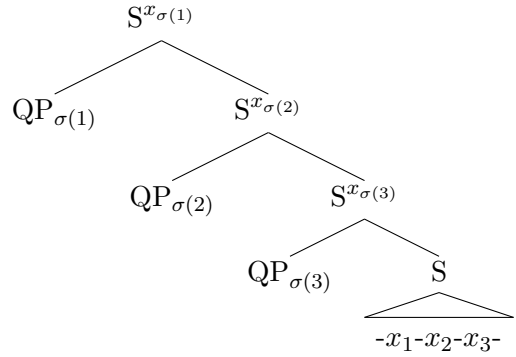
\Tree [.S$^{x_{\sigma(1)}}$ QP_{\sigma(1)} [.S$^{x_{\sigma(2)}}$ [.QP_{\sigma(2)} ][.S$^{x_{\sigma(3)}}$ [.QP_{\sigma(3)} ]\qroof{-$x_1$-$x_2$-$x_3$-}.S ]]]

\medskip
and the corresponding Computation Tree

\medskip
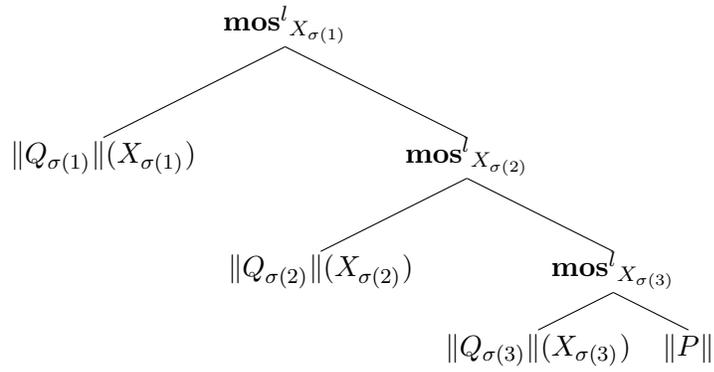
\Tree [.$\mosl_{X_{\sigma(1)}}$  [.$\|Q_{\sigma(1)}\|(X_{\sigma(1)})$ ][.$\mosl_{X_{\sigma(2)}}$  [.$\|Q_{\sigma(2)}\|(X_{\sigma(2)})$  ][.$\mosl_{X_{\sigma(3)}}$  [.$\|Q_{\sigma(3)}\|(X_{\sigma(3)})$  ] [.$\|P\|$  ]]]]
\par
\medskip
\end{flushleft}
The Computation Tree in (A3) gives rise to the following general map, with $\sigma\in S_3$ (where $S_3$ is the set of permutations of the set $\{1, 2, 3\}$)

\begin{center} \xext=1000 \yext=2500
\begin{picture}(\xext,\yext)(\xoff,\yoff)
   \put(-1000,2100){$\strat^{3,\sigma}_A:$}
  \putmorphism(500,2400)(0,-1)[\cC(X_1)\times\cP(X_1\times X_2\times X_3)\times\cC(X_2)\times\cC(X_3)` \phantom{\cC(X_1)\times\cC(X_2)\times\cC(X_3)\times\cP(X_1\times X_2\times X_3)}`\lk \bar{\pi}_{\sigma(1)},\bar{\pi}_{\sigma(2)}, \bar{\pi}_{\sigma(3)}, \pi_2 \rk ]{600}{1}r

  \putmorphism(500,1800)(0,-1)[\cC(X_{\sigma(1)})\times\cC(X_{\sigma(2)})\times\cC(X_{\sigma(3)})\times\cP(X_1\times X_2\times X_3)` \phantom{\cC(X_1)\times\cC(X_2)\times\cC\cP(X_1\times X_2\times X_3)}`1\times 1 \times \mosl_{X_{\sigma(3)}}]{600}{1}r

  \putmorphism(500,1200)(0,-1)[\cC(X_{\sigma(1)})\times\cC(X_{\sigma(2)})\times\cP(\ldots\times \widehat{X_{\sigma(3)}}\times \ldots)`\phantom{\cC(X)\times\cC(\cP(X))}`1\times \mosl_{X_{\sigma(2)}}]{600}{1}r

  \putmorphism(500,600)(0,-1)[\cC(X_{\sigma(1)})\times\cP(X_{\sigma(1)})`2`\mosl_{X_{\sigma(1)}}]{600}{1}r
 \end{picture}
\end{center}
There are six such maps corresponding to six permutations $\sigma$ of $\{ 1,2,3\}$. These maps are different in general.
Thus Strategy A yields 6 asymmetric readings for a sentence with three QPs.

\subsection{Strategy B}

In the polyadic approach (as developed in \cite{may85}, \cite{keenan87}, \cite{keenan92}, \cite{BZ}, \cite{benthem})
\begin{itemize}
\item Surface Structure Tree gets rewritten (disambiguated) as Formal Structure Trees (Polyadic Logical Forms) via
\begin{itemize}
  \item QR Rule;
  \item Predicate Collapsing;
  \item Rotation.
\end{itemize}
\item Formal Structure Trees (PLFs) are relabelled as Computation Trees as follows
\begin{itemize}
  \item $Polyadic$ (root of a subtree representing a polyadic quantifier) is interpreted as a suitably typed $\pu$-operation (we can choose globally whether we use only $\pul$ or $\pur$ and then consequently stick to it).
  \item S$^x$, S, QP are interpreted as above.
\end{itemize}
\end{itemize}

\vskip 2mm

\noindent \textbf{Sentence with one QP}, e.g. \textit{Every kid (most kids) entered.}

\begin{flushleft}
(B1) Surface Structure Tree

\medskip
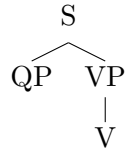
\Tree [.S [.QP ][.VP V ]]
\end{flushleft}

\begin{flushleft}
(B1) Formal Structure Tree (PLF) and the corresponding Computation Tree

\medskip
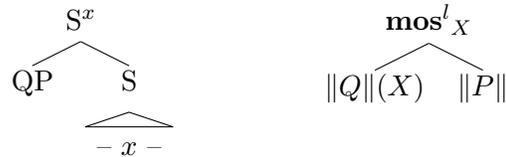
\Tree [.S$^x$ [.QP ] \qroof{-- $x$ --}.S ]
\Tree [.$\mosl_X$ $\|Q\|(X)$ $\|P\|$  ]
\par
\medskip
\end{flushleft}
The Computation Tree in (B1) gives rise to the following general map

\begin{center} \xext=600 \yext=600
\begin{picture}(\xext,\yext)(\xoff,\yoff)
   \put(-500,300){$\strat^1_B:$}
  \putmorphism(300,500)(0,-1)[\cC(X)\times\cP(X)`2`\mosl_X]{500}{1}r
 \end{picture}
\end{center}
In this case, there is one such map - thus Strategy B yields one reading for a sentence with one QP.

\vskip 2mm

\noindent \textbf{Sentence with two QPs}, e.g. \textit{Every girl likes a boy.}

\begin{flushleft}
(B2) Surface Structure Tree

\medskip
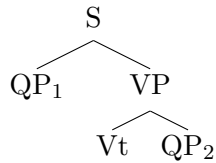
\Tree [.S [.QP_1 ][.VP [.Vt ][.QP_2 ]]]
\end{flushleft}

\medskip
\begin{flushleft}
(B2) Formal Structure Tree (PLF) obtained from LF in (A2) via rotation

\medskip
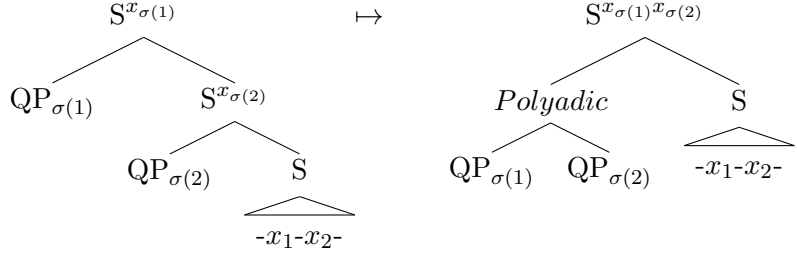
\Tree [.S$^{x_{\sigma(1)}}$ QP_{\sigma(1)} [.S$^{x_{\sigma(2)}}$ [.QP_{\sigma(2)} ] \qroof{-$x_1$-$x_2$-}.S ]]
$\mapsto$
\Tree [.S$^{x_{\sigma(1)} x_{\sigma(2)}}$ [.$Polyadic$ [.QP_{\sigma(1)} ][.QP_{\sigma(2)} ]] \qroof{-$x_1$-$x_2$-}.S ]
\par
\medskip

and the corresponding Computation Tree

\medskip

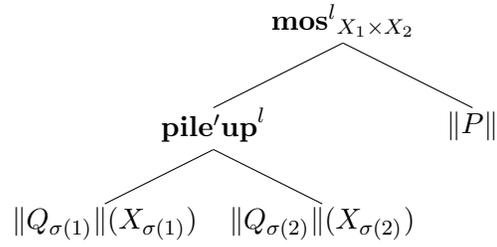
\Tree [.$\mosl_{X_{1}\times X_2}$ [.${\pul}$ [.$\|Q_{\sigma(1)}\|(X_{\sigma(1)})$  ][.$\|Q_{\sigma(2)}\|(X_{\sigma(2)})$  ]] $\|P\|$ ]

\par
\medskip
\end{flushleft}
The Computation Tree in (B2) gives rise to the following general map, with $\sigma\in S_2$

\medskip
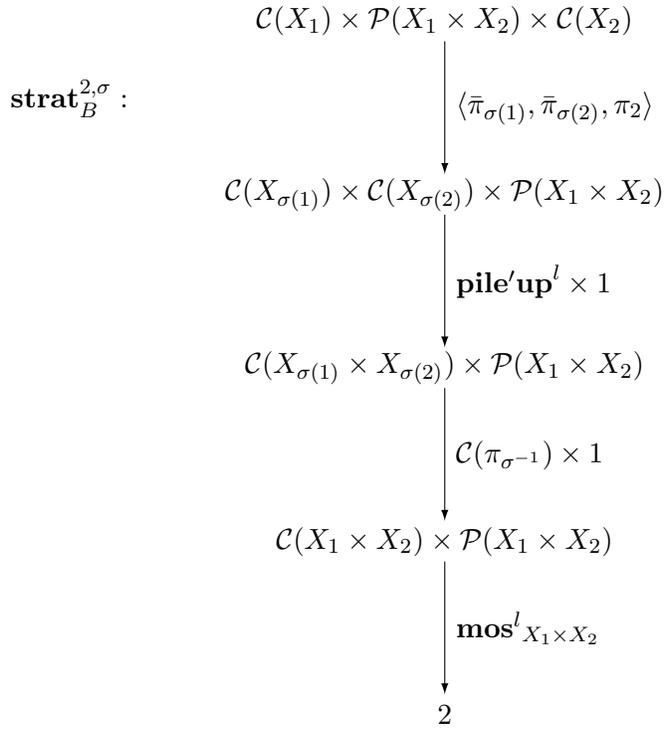
\begin{center} \xext=1000 \yext=2500
\begin{picture}(\xext,\yext)(\xoff,\yoff)
  \put(-1000,2100){$\strat^{2,\sigma}_B:$}
  \putmorphism(500,2400)(0,-1)[\cC(X_1)\times\cP(X_1\times X_2)\times\cC(X_2)`\phantom{\cC(X_1)\times\cC(X_2)\times\cP(X_1\times X_2)}`\lk \bar{\pi}_{\sigma(1)},\bar{\pi}_{\sigma(2)},\pi_2\rk]{600}{1}r

  \putmorphism(500,1800)(0,-1)[\cC(X_{\sigma(1)})\times\cC(X_{\sigma(2)})\times\cP(X_{1}\times X_{2})`\phantom{\cC(X)\times\cC(\cP(X)}`{\pul}\times 1 ]{600}{1}r

  \putmorphism(500,1200)(0,-1)[\cC(X_{\sigma(1)}\times X_{\sigma(2)})\times\cP(X_{1}\times X_{2})`\phantom{\cC(X)\times\cC(\cP(X))}`\cC(\pi_{\sigma^{-1}})\times 1 ]{600}{1}r
  \putmorphism(500,600)(0,-1)[\cC(X_1\times X_2)\times\cP(X_1\times X_2)`2`\mosl_{X_{1}\times X_2}]{600}{1}r
 \end{picture}
\end{center}

\vskip 2mm

\noindent There are two such maps corresponding to the two permutations $\sigma$ of $\{ 1,2\}$ combined with $\pul$-operation (in that case, alternatively, we can use $\pul$ and $\pur$). These maps are different in general. Thus Strategy B yields two (both) asymmetric readings for a sentence with two QPs.

\vskip 2mm

\noindent \textbf{Sentence with three QPs}, e.g. \textit{Some teacher gave every student most books.}
\begin{flushleft}
(B3) Surface Structure Tree

\medskip
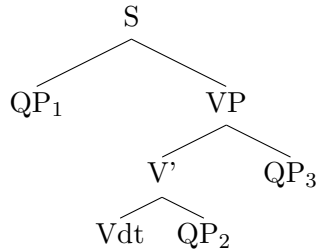
\Tree [.S [.QP_1 ][.VP [.V' [.Vdt ][.QP_2 ]][.QP_3 ]]]
\end{flushleft}

\begin{flushleft}
(B3) Formal Structure Tree (PLF) obtained from LF in (A3) via rotation

\medskip
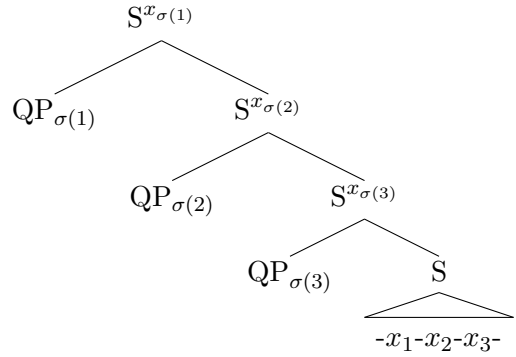
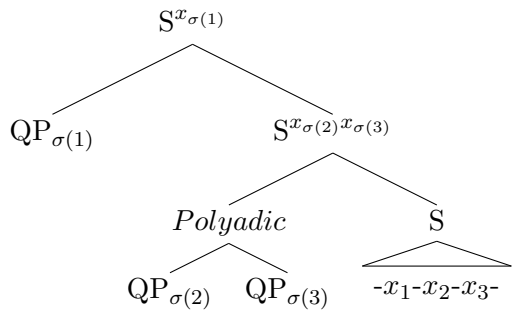
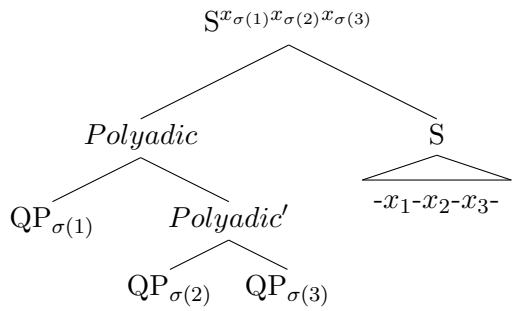
\Tree [.S$^{x_{\sigma(1)}}$ QP_{\sigma(1)} [.S$^{x_{\sigma(2)}}$ [.QP_{\sigma(2)} ][.S$^{x_{\sigma(3)}}$ [.QP_{\sigma(3)} ]\qroof{-$x_1$-$x_2$-$x_3$-}.S ]]]
\end{flushleft}

\begin{flushleft}
\medskip
$\mapsto$
\Tree [.S$^{x_{\sigma(1)}}$ QP_{\sigma(1)} [.S$^{x_{\sigma(2)} x_{\sigma(3)}}$ [.$Polyadic$ [.QP_{\sigma(2)} ][.QP_{\sigma(3)} ]] \qroof{-$x_1$-$x_2$-$x_3$-}.S ]]
\end{flushleft}

\begin{flushleft}
\medskip
$\mapsto$
\Tree [.S$^{x_{\sigma(1)} x_{\sigma(2)} x_{\sigma(3)}}$ [.$Polyadic$ QP_{\sigma(1)} [.$Polyadic'$ [.QP_{\sigma(2)} ][.QP_{\sigma(3)} ]]] \qroof{-$x_1$-$x_2$-$x_3$-}.S ]
\end{flushleft}

\begin{flushleft}
\medskip
and the corresponding Computation Tree
\medskip

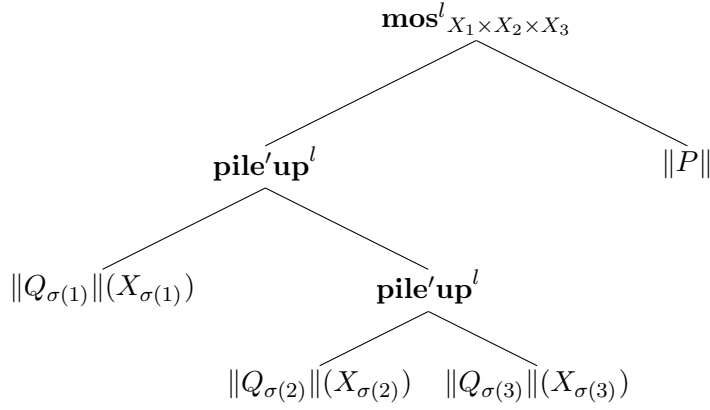
\Tree [.$\mosl_{X_{1}\times X_2\times X_3}$ [.${\pul}$  $\|Q_{\sigma(1)}\|(X_{\sigma(1)})$  [.${\pul}$ $\|Q_{\sigma(2)}\|(X_{\sigma(2)})$ $\|Q_{\sigma(3)}\|(X_{\sigma(3)})$ ]] $\|P\|$ ]
\par
\medskip
\end{flushleft}
The Computation Tree in (B3) gives rise to the following general map, with $\sigma\in S_3$

\begin{center} \xext=1000 \yext=3100
\begin{picture}(\xext,\yext)(\xoff,\yoff)
\put(-1000,2700){$\strat^{3,\sigma}_B:$}

  \putmorphism(500,3000)(0,-1)[\cC(X_1)\times\cP(X_1\times X_2\times X_3)\times\cC(X_2)\times\cC(X_3)` \phantom{\cC(X_1)\times\cC(X_2)\times\cC(X_3)\times\cC(\cP(X_1\times X_2\times X_3))}`\lk \bar{\pi}_{\sigma(1)},\bar{\pi}_{\sigma(2)},\bar{\pi}_{\sigma(3)},\pi_2\rk]{600}{1}r

  \putmorphism(500,2400)(0,-1)[\cC(X_{\sigma(1)})\times\cC(X_{\sigma(2)})\times\cC(X_{\sigma(3)})\times\cP(X_{1}\times X_{2}\times X_{3})`\phantom{\cC(X_1)\times\cC(X_2)\times\cP(X_1\times X_2\times X_3)}`1\times {\pul} \times 1]{600}{1}r

  \putmorphism(500,1800)(0,-1)[\cC(X_{\sigma(1)})\times\cC(X_{\sigma(2)}\times X_{\sigma(3)})\times\cP(X_1\times X_2\times X_3)`\phantom{\cC(X)\times\cC(\cP(X))}`{\pul} \times 1]{600}{1}r

   \putmorphism(500,1200)(0,-1)[\cC(X_{\sigma(1)}\times X_{\sigma(2)}\times X_{\sigma(3)})\times\cP(X_1\times X_2\times X_3)`\phantom{\cC(X)\times\cC(\cP(X))}`\cC(\pi_{\sigma^{-1}})\times 1]{600}{1}r

  \putmorphism(500,600)(0,-1)[\cC(X_1\times X_2\times X_3)\times\cP(X_1\times X_2\times X_3)`2`\mosl_{X_{1}\times X_2\times X_3}]{600}{1}r
 \end{picture}
\end{center}
There are six such maps corresponding to six permutations $\sigma$ of $\{ 1,2,3\}$ combined with $\pul$-operation (in that case we can choose globally whether we use only $\pul$ or $\pur$ and then consequently stick to it). These maps are different in general. Thus Strategy B yields 6 asymmetric readings for a sentence with three QPs.

\subsection{Strategy C}
In the continuation-based strategy approach (as proposed in \cite{barker})
\begin{itemize}
  \item Surface Structure Tree gets rewritten as Formal Structure Tree via
   \begin{itemize}
     \item no rewriting rules (Formal Structure Trees are just Surface Structure Trees - this is what is understood by in situ).
   \end{itemize}
  \item Relabelling Formal Structure Trees (= Surface Structure Trees) as the Computation Trees is as follows
\begin{itemize}
  \item S, VP, V' (roots of a (sub)tree with some (possibly all) arguments provided) are interpreted as suitably typed $\cps$-operations (left and right);
  \item V, Vt, Vdt (leafs of a tree) are interpreted as `continuized' (1-, 2-, 3-$ary$, respectively) predicates.
\end{itemize}
\end{itemize}

\vskip 2mm

\noindent \textbf{Sentence with one QP}, e.g. \textit{Every kid (most kids) entered.}

\begin{flushleft}
(C1) Surface Structure Tree and the corresponding Computation Tree

\medskip
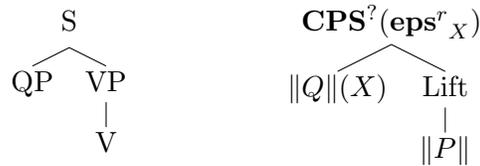
\Tree [.S [.QP ][.VP V ]]
\Tree [.$\cps^?(\epsr_{X})$ [.$\|Q\|(X)$  ][.Lift $\|P\|$ ]]
\par
\medskip
\end{flushleft}
The Computation Tree in (C1) gives rise to the following general map

\begin{center} \xext=600 \yext=1100
\begin{picture}(\xext,\yext)(\xoff,\yoff)
\put(-500,700){$\strat^1_c:$}
  \putmorphism(300,1000)(0,-1)[\cC(X)\times\cP(X)`\phantom{\cC(X)\times\cC\cP(X)}`1\times \eta_{\cP(X)}]{500}{1}r
  \putmorphism(300,500)(0,-1)[\cC(X)\times\cC\cP(X)`\cC(2)`\cps^?(\epsr_{X})]{500}{1}r

  \putmorphism(300,0)(1,0)[\phantom{\cC(2)}`2`ev_{id_2}]{600}{1}r
 \end{picture}
\end{center}
We use $\cps^?$ when it does not matter whether we apply $\cpsl$ or $\cpsr$. This is the case when one of the arguments is a lifted element (like interpretations of predicates in this strategy). Strategy C yields one reading for a sentence with one QPs.

\vskip 2mm

\noindent \textbf{Sentence with two QPs}, e.g. \textit{Every girl likes a boy.}

\begin{flushleft}
(C2) Surface Structure Tree and the corresponding Computation Tree

\medskip
\Tree [.S [.QP_1 ][.VP [.Vt ] [.QP_2 ]]]
\Tree [.$\cps^\varepsilon(\epsr_{X_1})$ [.$\|Q\|(X_1)$ ][.$\cps^?(\epsl_{X_2})$ [.Lift $\|P\|$ ][.$\|Q\|(X_2)$ ]]]
\par
\end{flushleft}
The Computation Tree in (C2) gives rise to the following general map

\begin{center} \xext=1000 \yext=1900
\begin{picture}(\xext,\yext)(\xoff,\yoff)
  \put(-1000,1500){$\strat^{2,\varepsilon}_C:$}
  \putmorphism(500,1800)(0,-1)[\cC(X_1)\times\cP(X_1\times X_2)\times\cC(X_2)`\phantom{\cC(X_1)\times\cC\cP(X_1\times X_2)\times\cC(X_2)}`1\times \eta_{\cP(X_1\times X_2)}\times 1]{600}{1}r
  \putmorphism(500,1200)(0,-1)[\cC(X_1)\times\cC\cP(X_1\times X_2)\times\cC(X_2)`\phantom{\cC(X)\times\cC\cP(X)}`1\times \cps^?(\epsl_{X_2})]{600}{1}r
  \putmorphism(500,600)(0,-1)[\cC(X_1)\times\cC\cP(X_1)`\cC(2)`\cps^\varepsilon(\epsr_{X_1})]{600}{1}r

  \putmorphism(500,0)(1,0)[\phantom{\cC(2)}`2`ev_{id_2}]{600}{1}r
 \end{picture}
\end{center}
with $\varepsilon\in \{l,r\}$. Depending on whether we use $\cpsl$ or $\cpsr$, we get either one or the other of the two asymmetric readings for a sentence with two QPs. Strategy C yields two readings for a sentence with two QPs corresponding to the two $\cps$'es.

\vskip 2mm

\noindent \textbf{Sentence with three QPs}, e.g. \textit{Some teacher gave every student most books.}

\begin{flushleft}
(C3) Surface Structure Tree and the corresponding Computation Tree

\medskip
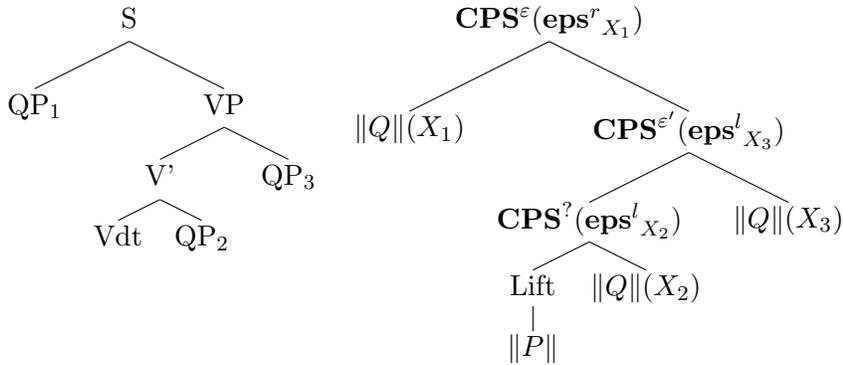
\Tree [.S [.QP_1 ][.VP [.V' [.Vdt ][.QP_2 ]][.QP_3 ]]]
\Tree [.$\cps^\varepsilon(\epsr_{X_1})$ [.$\|Q\|(X_1)$ ][.$\cps^{\varepsilon'}(\epsl_{X_3})$ [.$\cps^?(\epsl_{X_2})$ [.Lift $\|P\|$ ][.$\|Q\|(X_2)$ ]][.$\|Q\|(X_3)$ ]]]
\par
\end{flushleft}
The Computation Tree in (C3) gives rise to the following general map

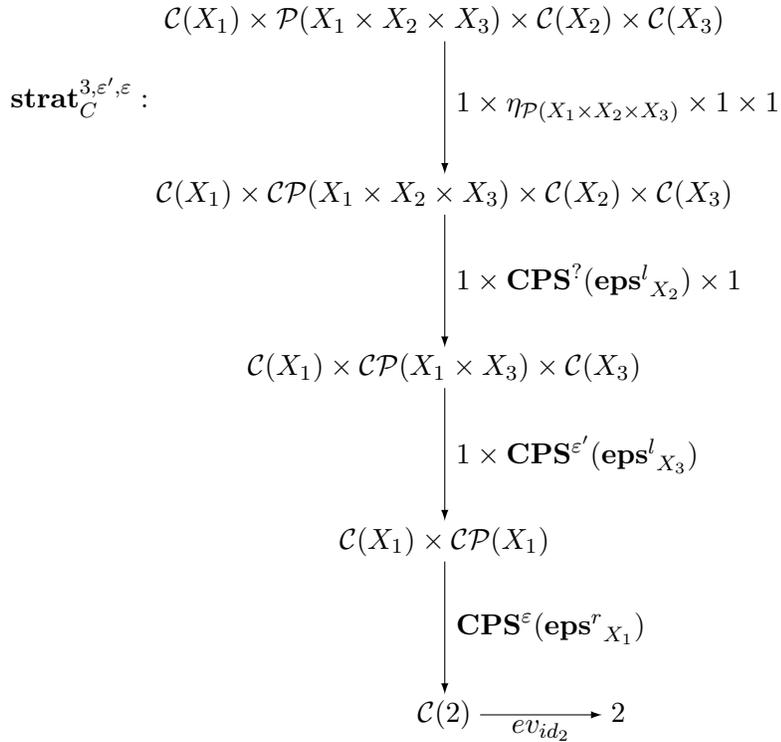
\begin{center} \xext=1000 \yext=2500
\begin{picture}(\xext,\yext)(\xoff,\yoff)
 \put(-1000,2100){$\strat^{3,\varepsilon',\varepsilon}_C:$}
  \putmorphism(500,2400)(0,-1)[\cC(X_1)\times\cP(X_1\times X_2\times X_3)\times\cC(X_2)\times\cC(X_3)` \phantom{\cC(X_1)\times\cC(X_2)\times\cC(X_3)\times\cC\cP(X_1\times X_2\times X_3)}`1\times\eta_{\cP(X_1\times X_2\times X_3)}\times 1\times 1]{600}{1}r

  \putmorphism(500,1800)(0,-1)[\cC(X_1)\times\cC\cP(X_1\times X_2\times X_3)\times\cC(X_2)\times\cC(X_3)`\phantom{\cC(X_1)\times\cC(X_3)\times\cC\cP(X_1\times X_2\times X_3)}`1\times \cps^?(\epsl_{X_2})\times 1]{600}{1}r

  \putmorphism(500,1200)(0,-1)[\cC(X_1)\times\cC\cP(X_1\times X_3)\times\cC(X_3)`\phantom{\cC(X)\times\cC\cP(X)}`1\times \cps^{\varepsilon'}(\epsl_{X_3})]{600}{1}r

  \putmorphism(500,600)(0,-1)[\cC(X_1)\times\cC\cP(X_1)`\cC(2)`\cps^\varepsilon(\epsr_{X_1})]{600}{1}r
  \putmorphism(500,0)(1,0)[\phantom{\cC(2)}`2`ev_{id_2}]{600}{1}r
 \end{picture}
\end{center}
Strategy C provides four asymmetric readings for the sentence such that QP in subject position can be placed either first or last only (corresponding to the four possible combinations of the two $\cps$'es). Thus it yields four out of six readings accounted for by strategies A and B.

The tables below summarize the main features of the three approaches.\\

\noindent \textbf{Passing from Surface Structure Tree Trees to Formal Structure Trees}

\[
\begin{array}{|c|c|c|c|} \hline

  \textit{\textbf{Strategy}} & \textit{\textbf{A}} & \textit{\textbf{B}}  & \textit{\textbf{C}} \\ \hline
  \textit{} & \textit{}  & \textit{} & \textit{} \\

  \textit{Rewrite} & \textit{QR,}  & \textit{QR,} & \textit{No rewrite rules} \\
  \textit{rules}&  \textit{Predicate} & \textit{Predicate} & \textit{(in situ)} \\
  \textit{} & \textit{Collapsing}  & \textit{Collapsing,} & \textit{} \\
  \textit{} & \textit{}  & \textit{Rotation} & \textit{} \\ \hline
 \end{array}
\]
\textbf{Passing from Formal Structure Trees to Computation Trees}
\[
\begin{array}{|c|c|c|c|} \hline
 \textit{\textbf{Strategy}} & \textit{\textbf{A}} & \textit{\textbf{B}}  & \textit{\textbf{C}} \\ \hline
 \textit{} & \textit{}  & \textit{} & \textit{} \\
  \textit{Relabelling} & \textit{$S^{x} \mapsto  \mos$} & \textit{$S^{\vec{x}} \mapsto  \mos$}  & \textit{} \\
  \textit{inner nodes} & \textit{}  & \textit{} & \textit{} \\
  \textit{} & \textit{} & \textit{$Polyadic \mapsto$}  & \textit{$S, V\! P, V' \mapsto$} \\
  \textit{} & \textit{} & \textit{$\pu$} & \textit{$\cps$ } \\ \hline

    \textit{} & \textit{}  & \textit{} & \textit{} \\
  \textit{Relabelling} & \textit{S $\mapsto$ relation} & \textit{S $\mapsto$ relation} & \textit{V, V\! t, V\! dt $\mapsto$}\\
   \textit{leaves} & \textit{}  & \textit{} & \textit{continuized relation} \\

  \textit{} & \textit{} & \textit{} & \textit{} \\
  \textit{ } & \textit{$QP\mapsto\; \cC$-comp.}  & \textit{$QP\mapsto\; \cC$-comp.} & \textit{$QP\mapsto\; \cC$-comp.} \\  \hline
\end{array}
\]
The semantics for sentences with intransitive and transitive verbs, as defined by the strategies A, B, and C, are equivalent. The semantics for sentences with ditransitive verbs, as defined by the strategies A, B, are equivalent. They provide six asymmetric readings of the sentence. The semantics for sentences with ditransitive verbs, as defined by the strategy C, provides four asymmetric readings of the sentence such that QP in subject position can be placed either first or last only. Thus they correspond to four out of six readings accounted for by strategies A and B. The proofs are given in the Appendix.

\section{Appendix}

\subsection{The continuation monad}
In this subsection, we gather all the basic facts (sometimes repeated from the text) of the {\em continuation monad} $\cC$ on $Set$. We have an adjunction
\begin{center}
\xext=1000 \yext=370
\begin{picture}(\xext,\yext)(\xoff,\yoff)
\putmorphism(0,140)(1,0)[Set`Set^{op}`]{1000}{0}a
\putmorphism(0,210)(1,0)[\phantom{Set}`\phantom{Set^{op}}`\cP]{1000}{1}a
\putmorphism(0,70)(1,0)[\phantom{Set}`\phantom{Set^{op}}`\cP^{op}]{1000}{-1}b
\end{picture}
\end{center}
where both $\cP$ and $\cP^{op}$ are the contravariant powerset functors\footnote{Note that this is in contrast with the functor $\cP$, where $\cP$ is the covariant power-set functor.} with the domains and codomains as displayed. In particular,
for $f:X\ra Y$, the function $\cP(f)=f^{-1} : \cP(Y)\ra \cP(X)$ is given by
  \[  f^{-1}(h) = h\circ f\]
for $h:Y\ra \bt$.\\
Function $\eta_X : X\ra \cC(X)$, the component at set $X$ of the unit of this adjunction $\eta : 1_{Set} \ra \cP\cP^{op}=\cC$, is given by
  \[ \eta_X(x)= \lambda h_{:\cP(X)}. h(x).\]
Function $\varepsilon_X : X\ra \cC(X)$, the component at set $X$ of the counit of this adjunction $\varepsilon : 1_{Set} \ra \cP^{op}\cP$, is given by (essentially the same formula)
  \[ \varepsilon_X(x)= \lambda h_{:\cP^{op}(X)}. h(x)\]
  for $x\in X$.\\
The function $\cC(f): \cC(X)\ra \cC(Y)$, for $Q:\cP(X)\ra \bt \in \cC(X)$, is a function $\cC(f)(Q):\cP(Y)\ra \bt$  given by
  \[ \cC(f)(Q)(h) = Q(h\circ f)\]
  for $h:Y\ra \bt$.\\

The monad induced by this adjunction is the continuation monad. Its multiplication is given by the counit of the above adjunction transported back to $Set$, i.e. $\mu = \cP^{op}(\varepsilon_{\cP})$. For $X$ in $Set$, the function
\[ \mu_X:\cC^2(X)\ra \cC(X)\]
is given by
 \[ \mu_X(\cR) = \cR\circ\eta_{\cP(X)}\]
  for $\cR\in \cC^{2}(X)$.\\
    In $\lambda$-notation we write
  \[\mu_X(\cF)(h)=\cF(\la D_{: \cC(X)}.D(h)).\]

The left strength for the monad $\cC$ is
 \[ \stl : \cC(X)\times Y \lra \cC(X\times Y) \]
 for $M\in \cC(X)$ and $y\in Y$, given by
 \[ \stl(M,y)=\la c_{:\cP(X\times Y)}. M(\la x_{:X}. c(x,y)): \cP(X\times Y)\ra \bt \]
  and the right strength, for $x\in X$ and $n\in \cC(Y)$, is given by
  \[ \str(x,N)=\la c_{: \cP(X\times Y)}. M(\lambda y_{:Y}.c(x,y)): \cP(X\times Y) \ra \bt.\]

The left pile'up operation
\[ \pul : \cC(X)\times \cC(Y) \lra \cC(X\times Y) \]
is the following composition
\begin{center} \xext=2700 \yext=250
\begin{picture}(\xext,\yext)(\xoff,\yoff)
  \putmorphism(0,50)(1,0)[\cC(X)\times\cC(Y)`\cC(X\times\cC(Y))`\stl]{900}{1}a
  \putmorphism(900,50)(1,0)[\phantom{\cC(X\times\cC(Y))}`\cC^2(X\times Y)`\cC(\str)]{900}{1}a
  \putmorphism(1800,50)(1,0)[\phantom{\cC^2(X\times Y)}`\cC(X\times Y)`\mu_{X\times Y}]{900}{1}a
 \end{picture}
\end{center}
where, for $Q\in \cC(X)$,  $Q'\in\cC(Y)$, $c\in \cP(X\times \cC(Y))$, we have
\[ \stl(Q,Q')(c)= Q(\la x_{:X} c(x,Q)) \]
and, for $d\in \cC(X\times \cC(Y))$, $\cU\in \cP\cC(X\times Y)$, we have
\[ \cC(\str)(d)(\cU) = d(\cU\circ \str).  \]

Now, using the above formulas, we can calculate $\pul$ as the composition on  $Q\in \cC(X)$,  $Q'\in\cC(Y)$, and $c\in \cP(X\times Y)$ as follows
\[ \pul(Q,Q')(c) = \]
\[ =\mu_{X\times Y}(\cC(\str)(\stl(Q,Q')))(c) = \]
\[ = \cC(\str)(\stl(Q,Q'))(\la D_{:\cC(X\times Y)} D(c)) = \]
\[ = \stl(Q,Q'))((\la D_{:\cC(X\times Y)} D(c))\circ \str) = \]
\[ = Q(\la x_{:X}((\la D_{:\cC(X\times Y)} D(c))\circ \str)(x,Q')) = \]
\[ = Q(\la x_{:X}((\la D_{:\cC(X\times Y)} D(c))(\str(x,Q')) = \]
\[ = Q(\la x_{:X} \str(x,Q')(c)) = \]
\[ = Q(\la x_{:X} Q'(\la y_{:Y}c(x,y))) \]
Similarly, we can show that
\[ \pur(Q,Q')(c) = Q'(\la y_{:Y} \Q(\la x_{:X}c(x,y))).\]

One can easily verify that pile'ups are related by
\[ \pur_{X,Y} = \cC(\pi_{(2,1)})\circ \pul_{Y,X}\circ \pi_{(2,1)}. \]

\subsection{Some properties of pile'up operations}

\begin{lemma}[Pile'up lemma]\label{app-pu-lem1}
$\pu$'s on pairs where one element is continuaized agree and are equal to the corresponding strength.
\end{lemma}

{\em Proof.} We have to show that the functions
\begin{center} \xext=1400 \yext=400
\begin{picture}(\xext,\yext)(\xoff,\yoff)
  \putmorphism(0,200)(1,0)[\phantom{T(X_1)\times  T(X_2)}`\phantom{T(X_1\times X_2)}`{\pul}_{X_1,X_2}]{1400}{1}a
  \putmorphism(0,150)(1,0)[T(X_1)\times  T(X_2)`T(X_1\times X_2)`]{1400}{0}a
  \putmorphism(0,100)(1,0)[\phantom{T(X_1)\times  T(X_2)}`\phantom{T(X_1\times X_2)}`{\pur}_{X_1,X_2}]{1400}{1}b
 \end{picture}
\end{center}
are equalized by both
\begin{center} \xext=1400 \yext=100
\begin{picture}(\xext,\yext)(\xoff,\yoff)
  \putmorphism(0,00)(1,0)[X_1\times  T(X_2)`T(X_1)\times T(X_2)`\eta_{X_1}\times T(1_{X_2})]{1400}{1}a
 \end{picture}
\end{center}
and
\begin{center} \xext=1400 \yext=100
\begin{picture}(\xext,\yext)(\xoff,\yoff)
  \putmorphism(0,00)(1,0)[T(X_1)\times  X_2`T(X_1)\times T(X_2)`T(1_{X_1})\times \eta_{X_2}]{1400}{1}a
 \end{picture}
\end{center}
and their composition with these functions are equal to strength morphisms.

Using the diagram
\begin{center} \xext=2600 \yext=2000
\begin{picture}(\xext,\yext)(\xoff,\yoff)
  \setsqparms[1`1`1`1;1300`600]
  \putsquare(0,1000)[T(X_1)\times X_2`T(X_1)\times T(X_2)`T(X_1\times X_2)`T(X_1\times T(X_2));T(1_{X_1})\times \eta_{X_2}`{\stl}_{X_1,X_2}`\stl_{T(X_1),X_2}`T(1_{X_1}\times \eta_{X_2})]
    \setsqparms[1`0`1`1;1300`600]
  \putsquare(1300,1000)[\phantom{T(X_1)\times T(X_2)}` T(T(X_1)\times X_2)`\phantom{T(X_1\times T(X_2))}`T^2(X_1\times X_2); \str_{T(X_1),X_2}``T({\stl}_{X_1,X_2})`T(\str_{X_1,X_2})]

  \putmorphism(2600,1000)(0,-1)[\phantom{T^2(X_1\times X_2)}`T(X_1\times X_2)`\mu_{X_1\times X_2}]{1000}{1}r

   \put(0,920){\line(0,-1){920}}
   \put(0,0){\vector(1,0){2300}}
    \put(1100,50){$1_{T(X_1\times X_2)}$}

    \put(200,920){\line(0,-1){320}}
    \put(200,600){\line(1,0){1800}}
   \put(2000,600){\vector(1,1){350}}
    \put(1100,650){$T(\eta_{X_1\times X_2})$}

    \put(100,920){\line(0,-1){520}}
    \put(100,400){\line(1,0){2100}}
   \put(2200,400){\vector(1,2){270}}
    \put(1100,310){$\eta_{T(X_1\times X_2)}$}

    \put(0,1700){\line(0,1){200}}
    \put(0,1900){\line(1,0){2600}}
   \put(2600,1900){\vector(0,-1){200}}
    \put(1100,1940){$\eta_{T(X_1)\times X_2}$}
 \end{picture}
\end{center}
we shall show that
\[ {\pur}_{X_1,X_2} \circ (T(1_{X_1})\times \eta_{X_2})= \stl_{X_1,X_2} =  {\pul}_{X_1,X_2}\circ (T(1_{X_1})\times \eta_{X_2}) .  \]
The other cases are symmetric. We have

\[ {\pur}_{X_1,X_2} \circ (T(1_{X_1})\times \eta_{X_2}) = \;\;({\rm def\; of}\;\pur) \]
\[ = \mu_{X_1,X_2}\circ T(\stl_{X_1,X_2})\circ \str_{T(X_1),X_2} \circ (T(1_{X_1})\times \eta_{X_2}) = \;\;(\eta\;{\rm strong\; w.r.t.}\; \str) \]
\[ = \mu_{X_1,X_2}\circ T(\stl_{X_1,X_2})\circ \eta_{T(X_1)\times X_2} = \;\; (\eta\; {\rm nat\; transf}) \]
\[ = \mu_{X_1,X_2}\circ \eta_{T(X_1\times X_2)})\circ \stl_{X_1,X_2} = \;\; (T\; {\rm monad}) \]
\[ = \stl_{X_1,X_2} \]
To show the remaining equation, we notice that we can continue the penultimate formula above as follows
\[{\pur}_{X_1,X_2} \circ (T(1_{X_1})\times \eta_{X_2}) =\ldots = \mu_{X_1,X_2}\circ \eta_{T(X_1\times X_2)})\circ \stl_{X_1,X_2} = \;\; (T\; {\rm monad}) \]
\[ = \mu_{X_1,X_2}\circ T(\eta_{X_1\times X_2})\circ \stl_{X_1,X_2} = \;\;(\eta\;{\rm strong\; w.r.t.}\; \str) \]
\[ = \mu_{X_1,X_2}\circ T(\str_{X_1,X_2})\circ   T(1_{X_1}\times \eta_{X_2})\circ \stl_{X_1,X_2} = \;\;(\stl\;{\rm nat\; transf}) \]
\[ = \mu_{X_1,X_2}\circ T(\str_{X_1,X_2})\circ \stl_{X_1,X_2}\circ   T(1_{X_1}\times \eta_{X_2}) = \;\;({\rm def\; of}\;\pul)  \]
\[ = {\pul}_{X_1,X_2} \circ (T(1_{X_1})\times \eta_{X_2})\]

$\diamondsuit$

\begin{corollary}
The left and right $\cps$-operation on pairs where one element is continuized agree.
\end{corollary}

{\em Proof.} The corollary states that, for any sets $X$, $Y$, $Z$ and a function $f:X\times Y\ra Z$, both morphisms
\begin{center} \xext=1200 \yext=100
\begin{picture}(\xext,\yext)(\xoff,\yoff)
  \putmorphism(0,00)(1,0)[X\times  T(Y)`T(X)\times T(Y)`\eta_{X}\times 1]{1200}{1}a
 \end{picture}
\end{center}
and
\begin{center} \xext=1200 \yext=100
\begin{picture}(\xext,\yext)(\xoff,\yoff)
  \putmorphism(0,00)(1,0)[T(X)\times  Y`T(X)\times T(Y)`1\times \eta_{Y}]{1200}{1}a
 \end{picture}
\end{center}
equalize the pair of morphisms
\begin{center} \xext=1400 \yext=400
\begin{picture}(\xext,\yext)(\xoff,\yoff)
  \putmorphism(0,200)(1,0)[\phantom{T(X)\times  T(Y)}`\phantom{T(Z)}`{\cpsl}(f)]{1400}{1}a
  \putmorphism(0,150)(1,0)[T(X)\times  T(Y)`Z`]{1400}{0}a
  \putmorphism(0,100)(1,0)[\phantom{T(X)\times  T(Y)}`\phantom{T(Z)}`{\cpsr}(f)]{1400}{1}b
 \end{picture}
\end{center}
This immediately follows from the above lemma and the definition of $\cps$'es. $\diamondsuit$


Using binary pile'up operations, we can define eight ternary pile'up operation
\[ T(X_1)\times  T(X_2)\times  T(X_3)\lra  T(X_1\times  X_2\times X_3) \]
out of the following diagram
\begin{center}
\xext=2600 \yext=1300 \adjust[`I;I`;I`;`I]
\begin{picture}(\xext,\yext)(\xoff,\yoff)
\put(600,1200){$T(X_1)\times T(X_2)\times T(X_3)$}
\put(0,600){$T^2(X_1\times X_2)\times T(X_3)$}
\put(1500,600){$T(X_1) \times T^2(X_2\times X_3)$}
\put(800,0){$T^3(X_1\times X_2\times X_3)$}
\put(200,500){\vector(2,-1){700}}
\put(400,500){\vector(2,-1){700}}
\put(2000,500){\vector(-2,-1){700}}
\put(2200,500){\vector(-2,-1){700}}

\put(150,350){$_{\pul}$}
\put(800,350){$_{\pur}$}

\put(1350,350){$_{\pul}$}
\put(2000,350){$_{\pur}$}

\put(1300,1100){\vector(2,-1){700}}
\put(1500,1100){\vector(2,-1){700}}
\put(900,1100){\vector(-2,-1){700}}
\put(1100,1100){\vector(-2,-1){700}}
\put(-50,850){$_{{\pul}\times 1}$}
\put(650,850){$_{{\pur}\times 1}$}
\put(1350,850){$_{1\times {\pul}}$}
\put(2050,850){$_{1\times {\pur}}$}
\end{picture}
\end{center}
However, both ${\pul}$ and ${\pur}$ operations are associative (Proposition \ref{app-pu-assoc} below) and hence only six of them are different, in general.

\begin{proposition}\label{app-pu-assoc}
Both $\pul$ and $\pur$ operations are associative on any monad on $Set$.
\end{proposition}

{\em Proof.} In fact, $\pul$ and $\pur$ are associative on any bi-strong monad on monoidal category.
We shall show this fact for a monad $T$ on $Set$ with the canonical strength.

We need to show that
\[ {\pur}\circ ({\pur}\times 1) = {\pur}\circ (1\times {\pur}) \]
and
\[ {\pul}\circ ({\pul}\times 1) = {\pul}\circ (1\times {\pul}) \]
But as pile'up are mutually definable, either of these equalities implies easily the other. We shall show the latter equality. For sets $X_1$, $X_2$, $X_3$, using all the assumptions, we have

\[ \pul_{X_1\times X_2,X_3}\circ (\pul_{X_1,X_2}\times 1_{T(X_3)}) =  \]

\[ = \mu_{X_1,\times X_2\times X_3} \circ  T(\str_{X_1\times X_2, X_3}) \circ  \underline{\stl_{X_1\times X_2, T(X_3)}\circ}\]
\[\underline{\circ\, (\mu_{X_1\times X_2}\times T(1_{X_3}))}\circ (T(\str_{X_1,X_2})\times 1_{T(X_3)})\circ (\stl_{X_1,T(X_2)}\times 1_{T(X_3)}) = \]

\[ = \mu_{X_1,\times X_2\times X_3} \circ  \underline{T(\str_{X_1\times X_2, X_3}) \circ  \mu_{X_1\times X_2\times T(X_3)}} \circ T(\stl_{X_1\times X_2, T(X_3)})\,\circ \]
\[ \circ\, \underline{\stl_{T(X_1\times X_2), T(X_3)}\circ (T(\str_{X_1,X_2})\times T(1_{X_3}))}\circ (\stl_{X_1,T(X_2)}\times 1_{T(X_3)}) = \]

\[ = \mu_{X_1,\times X_2\times X_3} \circ  \mu_{T(X_1\times X_2\times X_3)}  \circ  T^2(\str_{X_1\times X_2, T(X_3)})\circ T(\stl_{X_1\times X_2, T(X_3)})\,\circ \]
\[ \circ\, (T(\str_{X_1,X_2}\times 1_{X_3})) \circ\, \underline{\stl_{T(X_1\times X_2), T(X_3)}\circ (\stl_{X_1,T(X_2)}\times 1_{T(X_3)})} = \]

\[ = \mu_{X_1,\times X_2\times X_3} \circ  \mu_{T(X_1\times X_2\times X_3)}  \circ  T^2(\str_{X_1\times X_2, T(X_3)})
\circ \underline{T(\stl_{X_1\times X_2, T(X_3)})\,\circ }\]
\[ \underline{\circ\, (T(\str_{X_1,X_2}\times 1_{T(X_3)}))} \circ\, \stl_{X_1, T(X_2)\times T(X_3)} = \]

\[ = \mu_{X_1,\times X_2\times X_3} \circ  \mu_{T(X_1\times X_2\times X_3)}  \circ  T^2(\str_{X_1\times X_2, T(X_3)})\circ T(\str_{X_1, X_2\times T(X_3)})\,\circ \]
\[ \circ\, \underline{T(1_{X_1}\times \stl_{X_2, T(X_3)}) \circ\, \stl_{X_1, T(X_2)\times T(X_3)}} = \]

\[ = \mu_{X_1,\times X_2\times X_3} \circ  \mu_{T(X_1\times X_2\times X_3)}  \circ  \underline{T^2(\str_{X_1\times X_2, T(X_3)})}\circ T(\str_{X_1,X_2\times T(X_3)})\,\circ \]
\[ \circ\, \stl_{X_1, T(X_2\times T(X_3))} \circ (T(1_{X_1})\times \stl_{X_2,T(X_3)})  = \]

\[ = \mu_{X_1,\times X_2\times X_3} \circ  \mu_{T(X_1\times X_2\times X_3)}  \circ  T^2(\str_{X_1,T(X_2\times X_3)})\circ\underline{ T^2(1_{X_1}\times \str_{X_2,X_3})\,\circ} \]
\[ \underline{\circ\, T(\str_{X_1,X_2\times T(X_3)})} \circ\, \stl_{X_1, T(X_2\times T(X_3))} \circ (T(1_{X_1})\times \stl_{X_2,T(X_3)})  = \]

\[ = \mu_{X_1,\times X_2\times X_3} \circ  \mu_{T(X_1\times X_2\times X_3)}  \circ  T^2(\str_{X_1,T(X_2\times X_3)})
\circ T(\str_{X_1,T(X_2\times X_3)})\, \circ \]
\[ \circ\, \underline{T(1_{X_1}\times T(\str_{X_2,X_3}))\circ \stl_{X_1, T(X_2\times T(X_3))}} \circ (T(1_{X_1})\times \stl_{X_2,T(X_3)})  = \]

\[ = \mu_{X_1,\times X_2\times X_3} \circ \underline{ \mu_{T(X_1\times X_2\times X_3)}  \circ  T^2(\str_{X_1,T(X_2\times X_3)})
\circ T(\str_{X_1,T(X_2\times X_3)})} \;\circ\, \]
\[ \circ\, \stl_{X_1, T^2(X_2\times X_3)} \circ\, (T(1_{X_1})\times T(\str_{X_2,X_3})) \circ (T(1_{X_1})\times \stl_{X_2,T(X_3)})  = \]

\[ = \mu_{X_1,\times X_2\times X_3} \,\circ\, T(\str_{X_1,X_2\times X_3}) \circ\,\underline{ T(1_{X_1}\times \mu_{ X_2\times X_3} ) \;\circ} \]
\[ \underline{\circ\, \stl_{X_1, T^2(X_2\times X_3)}} \circ\, (T(1_{X_1})\times T(\str_{X_2,X_3})) \circ (T(1_{X_1})\times \stl_{X_2,T(X_3)})  = \]

\[ = \mu_{X_1,\times X_2\times X_3} \,\circ\, T(\str_{X_1,X_2\times X_3}) \circ\, \stl_{X_1, T(X_2\times X_3)}\circ\ \]
\[ \circ\, (T(1_{X_1})\times \mu_{ X_2\times X_3} )  \circ\, (T(1_{X_1})\times T(\str_{X_2,X_3})) \circ (T(1_{X_1})\times \stl_{X_2,T(X_3)})  = \]

\[ =\pul_{X_1,X_2\times X_3}\circ (1_{T(X_3)}\times \pul_{X_2,X_3}) \]
$\diamondsuit$

\subsection{Arity one: intransitive verbs}

\begin{proposition}
The semantics for sentences with intransitive verbs, as defined by the strategies A, B, and C, are equivalent.
\end{proposition}

{\em Proof.} In case of a sentence with an intransitive verb the semantics are defined by morphisms  $\strat^1_A$, $\strat^1_B$, and $\strat^1_C$. We need to show that they are equal. We have
\[ \strat^1_A= \mosl_X= \strat^1_B. \]

$\strat^1_C$ is the composition of the following morphisms
\begin{center} \xext=2600 \yext=200
\begin{picture}(\xext,\yext)(\xoff,\yoff)
  \putmorphism(0,50)(1,0)[\cC(X)\times\cP(X)`\cC(X)\times\cC\cP(X)`1\times \eta_{\cP(X)}]{1200}{1}a
  \putmorphism(1200,50)(1,0)[\phantom{\cC(X)\times\cC\cP(X)}`\cC(\bt)`\cpsl(\epsr_X)]{1200}{1}a
  \putmorphism(2400,50)(1,0)[\phantom{\cC(\bt)}`\bt`ev_{id_\bt}]{500}{1}a
 \end{picture}
\end{center}
Thus we need to show that this composition is equal to $\mosl_X$. Consider the following diagram
\begin{center}
\xext=2000 \yext=800
\begin{picture}(\xext,\yext)(\xoff,\yoff)
\setsqparms[1`1`-1`0;2000`600]
\putsquare(0,100)[\cC(X)\times \cP(X)`\bt`\cC(X)\times \cC\cP(X)`\cC(\bT);\mosl_X`1\times\eta_{\cP(X)}`ev_{id_\bt}`]

\putmorphism(0,100)(1,0)[\phantom{\cC(X)\times \cC\cP(X)}`\cC(X\times \cP(X))`\pul_X]{1000}{1}b
\putmorphism(1000,100)(1,0)[\phantom{\cC(X\times \cP(X))}`\phantom{\cC(\bt)}`\cC(\epsr_X]{1000}{1}b
\putmorphism(0,650)(2,-1)[\phantom{\cC(X)\times \cP(X)}`\phantom{\cC(X\times \cP(X))}`\stl]{1000}{1}l
\putmorphism(1200,300)(2,1)[\phantom{\cC(X\times \cP(X))}`\phantom{\bt}` ev_{\epsr_X}]{500}{1}r

\end{picture}
\end{center}

The left triangle commutes, as a consequence of Lemma \ref{app-pu-lem1}. To see that the mid triangle commutes, we take $M\in \cC(X)$ and $h\in \cP(X)$, and calculate
\[ ev_{\epsr_X}\circ \str(Q,h) = \]
\[ = ev_{\epsr_X}(\la D_{:\cP(X\times\cP(X))}M(\la x_{:X} D(x,h))) = \]
\[ = M(\la x_{:X} \epsr_X(x,h)) = \]
\[ = M(\la x_{:X} h(x)) =\]
\[ = N(h)=\mosl(N,h). \]

Finally, to see that the right triangle commutes, we take $N\in \cC(X\times \cP(X))$ and calculate
\[ ev_{id_\bt}\circ\cC(\epsr_X)(N)= \]
\[ = ev_{id_\bt} (\la c_{:\cP(\bt)} N(c\circ \epsr_X)) = \]
\[ = N( \epsr_X) = ev_{\epsr_X}(N). \]
Thus the whole diagram commutes, and hence $\strat^1_C=\mosl_X$, as required.
$\diamondsuit$

The above proof shows also the following technical lemma.

\begin{lemma}\label{app-mos-st}
For any set $X$, the diagram
\begin{center}
\xext=1000 \yext=600
\begin{picture}(\xext,\yext)(\xoff,\yoff)
\setsqparms[1`1`-1`1;1000`450]
\putsquare(0,50)[\cC(X)\times \cP(X)`\bt`\cC(X\times \cP(X))`\cC(\bt);\mosl_X`\stl`ev_{id_\bt}`\cC(\epsr_X)]
\end{picture}
\end{center}
commutes.
\end{lemma}

\subsection{Arity two: transitive verbs}

\begin{proposition}
The semantics for sentences with transitive verbs, as defined by the strategies A, B, and C, are equivalent. They provide two asymmetric readings of the sentence.
\end{proposition}

{\em Proof.}  In case of sentences with transitive verbs the semantics are defined by morphisms  $\strat^{2,\sigma}_A$, $\strat^{2,\sigma}_B$, and $\strat^{2,\varepsilon}_C$, with $\sigma\in S_2=\{id_2,\tau\}$ and $\varepsilon\in \{ l,r\}$. We need to show the equalities
\[  \strat^{2,\sigma}_A = \strat^{2,\sigma}_B,\]
for $\sigma\in S_2$, and
\[  \strat^{2,id_2}_B =\strat^{2,l}_C,\hskip 5mm     \strat^{2,\tau}_B = \strat^{2,r}_C.\]

To show the first equality, with $Q_1\in \cC(X_1)$, $Q_2\in \cC(X_2)$, and $P\in \cP(X_1\times X_2)$, we have
\[ \strat^{2,\sigma}_A(Q_1,Q_2,P) =\]
\[ = \mosl_{X_{\sigma(1)}}(Q_{\sigma(1)},\mosl_{X_{\sigma(2)}}(Q_{\sigma(2)},P)) =\]
\[ = \mosl_{X_{\sigma(1)}}(Q_{\sigma(1)},\la x_{\sigma(1): X_{\sigma(1)}}. Q_{\sigma(2)}(\la x_{\sigma(2):X_{\sigma(2)}}.P(x_1,x_2))) =\]
\[ = Q_{\sigma(1)}(\la x_{\sigma(1): X_{\sigma(1)}}. Q_{\sigma(2)}(\la x_{\sigma(2):X_{\sigma(2)}}.P(x_1,x_2))) =\]
\[ = Q_{\sigma(1)}(\la x_{\sigma(1): X_{\sigma(1)}}. Q_{\sigma(2)}(\la x_{\sigma(2):X_{\sigma(2)}}.P(\pi_{\sigma^{-1}}(x_{\sigma(1)},x_{\sigma(2)}))) =\]
\[ = \pul(Q_{\sigma(1)},Q_{\sigma(2)})(P\circ\pi^{-1}) = \]
\[ \cC(\pi_{\sigma^{-1}})(\pul(Q_{\sigma(1)},Q_{\sigma(2)}))(P) =\]
\[ =\strat^{2,\sigma}_B(Q_1,Q_2,P) \]

To show the remaining two equalities, let us first note that if either $\sigma=id_2$ and $\varepsilon=l$ or $\sigma=\tau$ and $\varepsilon=r$, we have
\[ \pu^\varepsilon= \cC(\pi_{\sigma^{-1}})\circ \pul\circ \pi_\sigma. \]
Thus we shall assume the above equation relating $\sigma$ with $\varepsilon$, and, with $Q_1\in \cC(X_1)$, $Q_2\in \cC(X_2)$, and $P\in \cP(X_1\times X_2)$, we obtain (the diagram illustrating these calculations would be too big to fit a page but the reader is encouraged to draw one)

\[  \strat^{2,\varepsilon}_C = \]
\[ = ev_{id_\bt} \circ \cps^\varepsilon(\epsr_{X_1})\circ (1\times\cps^?(\epsr_{X_2}))\circ (1\times 1\times \eta_{\cP(X_1\times X_2)} = \]
\[ = ev_{id_\bt} \circ  \cC(\epsr_{X_1})\circ \pu^\varepsilon\circ (\cC(1)\times\cC(\epsr_{X_2}))\circ (1\times \pu^?)\circ (1\times 1\times \eta) = \]
\[ = ev_{id_\bt} \circ  \cC(\epsr_{X_1})\circ \cC(1\times\epsr_{X_2})\circ \pu^\varepsilon\circ (1\times \pu^?)\circ (1\times 1\times \eta) = \]
\[ = ev_{id_\bt} \circ  \cC(\epsr_{X_1\times X_2})\circ \pu^\varepsilon\circ (1\times \pu^?)\circ (1\times 1\times \eta) = \]
\[ = ev_{id_\bt} \circ  \cC(\epsr_{X_1\times X_2})\circ \pu^?\circ (\pu^\varepsilon\times 1)\circ (1\times 1\times \eta) = \]
\[ = ev_{id_\bt} \circ  \cC(\epsr_{X_1\times X_2})\circ \pu^?\circ (1\times \eta)\circ (1\times \pu^\varepsilon) = \]
\[ = ev_{id_\bt} \circ  \cC(\epsr_{X_1\times X_2})\circ \stl\circ (1\times \pu^\varepsilon) = \]
\[ = ev_{id_\bt} \circ  \cC(\epsr_{X_1\times X_2})\circ \stl\circ (\cC(\pi_{\sigma^{-1}})\times 1)\circ (\pul\times 1)\circ (\pi_{\sigma}\times 1) = \]
\[ = \mosl_{X_1\times X_2}\circ (\cC(\pi_{\sigma^{-1}})\times 1)\circ (\pul\times 1)\circ (\pi_{\sigma}\times 1) = \]
\[ = \strat^{2,\sigma}_B \]

In the above calculations we used: the definition of $\cps$'es, naturality of $\pu^\varepsilon$, relations between $\eps$ morphisms, associativity of $\pu^\varepsilon$ (Proposition \ref{app-pu-assoc}), properties of product morphisms, pile'up lemma, and finally Lemma \ref{app-mos-st}.

Here and below $\cps^?$, $\pu^?$ stands for either $\cps^l$, $\pu^l$ or $\cps^r$, $\pu^r$ whatever is more convenient at the moment as it does not influence the end result. $\diamondsuit$

\subsection{Arity three: ditransitive verbs}

\begin{proposition}
The semantics for sentences with ditransitive verbs, as defined by the strategies A, B, are equivalent. They provide six asymmetric readings of the sentence.
\end{proposition}

{\em Proof.} In case of sentences with ditransitive verbs the semantics are defined by morphisms  $\strat^{3,\sigma}_A$, $\strat^{3,\sigma}_B$, and $\strat^{2,\varepsilon}_C$, with $\sigma\in S_3$ and $\varepsilon\in \{ l,r\}$. We need to show the equalities
\[  \strat^{3,\sigma}_A = \strat^{3,\sigma}_B,\]
for $\sigma\in S_3$.

The calculations are similar to those for transitive verbs. We present them for completeness. With $Q_1\in \cC(X_1)$, $Q_2\in \cC(X_2)$, $Q_3\in \cC(X_3)$, and $P\in \cP(X_1\times X_2\times X_3)$, we have
\[ \strat^{3,\sigma}_A(Q_1,Q_2,Q_3,P) =\]
\[ = \mosl_{X_{\sigma(1)}}(Q_{\sigma(1)},\mosl_{X_{\sigma(2)}}(Q_{\sigma(2)},\mosl_{X_{\sigma(3)}}(Q_{\sigma(3)},P)) =\]
\[ = Q_{\sigma(1)}(\la x_{\sigma(1): X_{\sigma(1)}}. Q_{\sigma(2)}(\la x_{\sigma(2):X_{\sigma(2)}}.
Q_{\sigma(3)}(\la x_{\sigma(3): X_{\sigma(3)}}.P(x_1,x_2,x_3))) =\]
\[ = Q_{\sigma(1)}(\la x_{\sigma(1): X_{\sigma(1)}}. Q_{\sigma(2)}(\la x_{\sigma(2):X_{\sigma(2)}}. Q_{\sigma(3)}(\la x_{\sigma(3): X_{\sigma(3)}}. P(\pi_{\sigma^{-1}}(x_{\sigma(1)},x_{\sigma(2)},x_{\sigma(3)}))) =\]
\[ = \pul(Q_{\sigma(1)},\pul(Q_{\sigma(2)},Q_{\sigma(3)}))(P\circ\pi_{\sigma^{-1}}) = \]
\[ \cC(\pi_{\sigma^{-1}})(\pul(Q_{\sigma(1)},\pul(Q_{\sigma(2)},Q_{\sigma(3)})))(P) =\]
\[ =\strat^{2,\sigma}_B(Q_1,Q_2,Q_3,P) \]
as required. $\diamondsuit$

\begin{proposition}
The semantics for sentences with ditransitive verbs, as defined by the strategy C, provides four asymmetric readings of the sentence such that QP in subject position can be placed either first or last only. Thus they correspond to four out of six readings accounted for by strategies A and B.
\end{proposition}

{\em Proof.}  In case of sentences with ditransitive verbs the semantics, according to strategies B and C, are defined by morphisms $\strat^{3,\sigma}_B$, $\strat^{3,\varepsilon,\varepsilon'}_C$, respectively. As we shall show, these morphisms are equal whenever $\sigma\in S_3$ is related to the pair $\lk\varepsilon', \varepsilon\rk\in \{l,r\}^2$ via relation
\[ \pu^{\varepsilon'}\circ (1\times \pu^{\varepsilon})= \cC(\pi_{\sigma^{-1}}) \circ \pul\circ (1\times \pul) \circ \pi_{\sigma} \]

As $\pul$ leaves the order intact and $\pur$ swaps the order, we can see that we have the following correspondence

\[
  \begin{array}{|c|c|} \hline
    \sigma & \lk \varepsilon',\varepsilon\rk \\ \hline
    (1,2,3) & \lk l,l\rk \\
    (1,3,2) & \lk l,r\rk \\
     (2,3,1)& \lk r,l\rk \\
     (3,2,1)& \lk r,r\rk \\
     (2,1,3)& - \\
     (3,1,2)&  -\\ \hline
  \end{array}
\]
Thus we shall assume the $\sigma$ is related to the pair $\lk \varepsilon,\varepsilon'\rk$, and, with $Q_1\in \cC(X_1)$, $Q_2\in \cC(X_2)$, $Q_3\in \cC(X_3)$, and $P\in \cP(X_1\times X_2\times X_3)$, we obtain (the diagram illustrating these calculations would be again too big to fit a page but the reader is encouraged to draw one)

\[  \strat^{3,\varepsilon',\varepsilon}_C = \]

\[ = ev_{id_\bt} \circ \cps^{\varepsilon'}(\epsr_{X_1})\circ  (1\times\cps^{\varepsilon}(\epsr_{X_2}))\circ\]
\[ \circ (1\times 1\times\cps^?(\epsr_{X_3}))\circ (1\times 1\times 1\times \eta) = \]

\[ = ev_{id_\bt} \circ \cC(\epsr_{X_1})\circ \pu^{\varepsilon'}\circ  (\cC(1)\times\cC(\epsr_{X_2}))\circ (1\times \pu^\varepsilon)\circ \]
\[ \circ (\cC(1)\times \cC(1)\times \cC(\epsr_{X_3}))\circ (1\times 1\times \pu^?)\circ (1\times 1\times 1\times \eta) = \]

\[ = ev_{id_\bt} \circ \cC(\epsr_{X_1})\circ  (\cC(1\times\epsr_{X_2}))\circ \pu^{\varepsilon'}\circ (\cC(1)\times \cC(1\times \epsr_{X_3}))\circ  \]
\[ \circ(1\times \pu^\varepsilon)\circ (1\times 1\times \pu^?)\circ (1\times 1\times 1\times \eta) = \]

\[ = ev_{id_\bt} \circ \cC(\epsr_{X_1})\circ  (\cC(1\times\epsr_{X_2}))\circ (\cC(1\times 1\times \epsr_{X_3}))\circ \pu^{\varepsilon'}\circ  \]
\[ \circ(1\times \pu^\varepsilon)\circ (1\times 1\times \pu^?)\circ (1\times 1\times 1\times \eta) = \]

\[ = ev_{id_\bt} \circ \cC(\epsr_{X_1\times X_2\times X_3})\circ \pu^{\varepsilon'}\circ  \]
\[ \circ(1\times \pu^\varepsilon)\circ (1\times 1\times \pu^?)\circ (1\times 1\times 1\times \eta) = \]

\[ = ev_{id_\bt} \circ \cC(\epsr_{X_1\times X_2\times X_3})\circ \pu^{\varepsilon'}\circ  \]
\[ \circ(1\times \pu^\varepsilon)\circ (1\times 1\times \pu^?)\circ (1\times 1\times 1\times \eta) = \]

\[ = ev_{id_\bt} \circ \cC(\epsr_{X_1\times X_2\times X_3})\circ \pu^{\varepsilon'}\circ  \]
\[ \circ (1\times \pu^?)\circ(1\times \pu^\varepsilon\times 1)\circ\ (1\times 1\times 1\times \eta) = \]

\[ = ev_{id_\bt} \circ \cC(\epsr_{X_1\times X_2\times X_3})\circ \pu^?\circ  \]
\[ \circ (\pu^{\varepsilon'}\times 1)\circ(1\times \pu^\varepsilon\times 1)\circ\ (1\times 1\times 1\times \eta)\stackrel{*}{=} \]

\[ \stackrel{*}{=} ev_{id_\bt} \circ \cC(\epsr_{X_1\times X_2\times X_3})\circ \pu^?\circ  \]
\[ \circ (\cC(\pi_{\sigma^{-1}})\times \cC(1)) \circ (\pul\times 1)\circ (1\times \pul \times 1) \circ (\pi_{\sigma}\times 1)\circ\ (1\times 1\times 1\times \eta) = \]

\[ = ev_{id_\bt} \circ \cC(\epsr_{X_1\times X_2\times X_3})\circ (\cC(\pi_{\sigma^{-1}}\times 1)\circ \pu^?\circ  \]
\[  \circ (\pul\times 1)\circ (1\times \pul \times 1) \circ (\pi_{\sigma}\times 1)\circ\ (1\times 1\times 1\times \eta) = \]

\[ = ev_{id_\bt} \circ \cC(\epsr_{X_1\times X_2\times X_3})\circ (\cC(\pi_{\sigma^{-1}}\times 1)\circ \pu^?\circ  \]
\[  \circ (\pul\times 1)\circ (1\times \pul \times 1) \circ\ (1\times 1\times 1\times \eta)\circ (\pi_{\sigma}\times 1) = \]

\[ = ev_{id_\bt} \circ \cC(\epsr_{X_1\times X_2\times X_3})\circ (\cC(\pi_{\sigma^{-1}}\times 1)\circ \pu^?\circ  \]
\[  \circ (\pul\times 1) \circ\ (1\times \times \eta)\circ (1\times \pul \times 1)\circ (\pi_{\sigma}\times 1) = \]

\[ = ev_{id_\bt} \circ \cC(\epsr_{X_1\times X_2\times X_3})\circ (\cC(\pi_{\sigma^{-1}}\times 1)\circ \pu^?\circ  \]
\[ = \circ (\pul\times 1) \circ\ (1\times 1\times \eta)\circ (1\times \pul \times 1)\circ (\pi_{\sigma}\times 1) = \]

\[ = ev_{id_\bt} \circ \cC(\epsr_{X_1\times X_2\times X_3})\circ (\cC(\pi_{\sigma^{-1}}\times 1)\circ \pu^?\circ  \]
\[   \circ\ (1\times \eta)\circ (\pul\times 1) \circ (1\times \pul \times 1)\circ (\pi_{\sigma}\times 1) = \]

\[ = ev_{id_\bt} \circ \cC(\epsr_{X_1\times X_2\times X_3})\circ (\cC(\pi_{\sigma^{-1}}\times 1)\circ \stl^l \circ  \]
\[   \circ (\pul\times 1) \circ (1\times \pul \times 1)\circ (\pi_{\sigma}\times 1) = \]

\[ = ev_{id_\bt} \circ \cC(\epsr_{X_1\times X_2\times X_3})\circ \stl^l \circ(\cC(\pi_{\sigma^{-1}})\times \cC(1))\circ   \]
\[   \circ (\pul\times 1) \circ (1\times \pul \times 1)\circ (\pi_{\sigma}\times 1) = \]

\[ = \mos^l_{X_1\times X_2\times X_3} \circ(\cC(\pi_{\sigma^{-1}})\times \cC(1))\circ   \]
\[   \circ (\pul\times 1) \circ (1\times \pul \times 1)\circ (\pi_{\sigma}\times 1) = \]

\[  = \strat^{3,\sigma}_B \]

In the above calculations we used: the definition of $\cps$'es, naturality of $\pu$'s (four times in three non-consecutive steps!), relations between $\eps$ morphisms, associativity of $\pu$'s (Proposition \ref{app-pu-assoc}), relations between $\sigma$ and $\lk \varepsilon',\varepsilon\rk$, properties of product morphisms (three consecutive steps), pile'up lemma, naturality of strength,  and finally Lemma \ref{app-mos-st}. $\diamondsuit$

\end{document}